\documentclass{article}

\title{Lifts of longest elements to braid groups acting on derived categories}
\author{Joseph Grant\\School of Mathematics, Leeds University, UK\\\texttt{joseph.grant@gmail.com}}

\usepackage{latexsym}
\usepackage{verbatim}
\usepackage{amsmath}
\usepackage{amssymb}
\usepackage{mathrsfs}
\usepackage{amsthm}
\usepackage{sfmath}

\input xy
\xyoption{all}
\newcommand{\shift}{\Sigma}

\newcommand{\be}{\begin{equation}}
\newcommand{\ee}{\end{equation}}

\newcommand{\cma}{\text{,}}
\newcommand{\fs}{\text{.}}
\newcommand{\semic}{\text{;}}

\newcommand{\da}{\text{-}}

\newcommand{\Hom}{\operatorname{Hom}\nolimits}

\newcommand{\dert}{\otimes^\textbf{L}}

\newcommand{\Aut}{\operatorname{Aut}\nolimits}

\newcommand{\per}{\operatorname{per}\nolimits}
\newcommand{\cone}{\operatorname{cone}\nolimits}

\newcommand{\fa}{\text{ for all }}
\newcommand{\ev}{\operatorname{ev}\nolimits}

\newcommand{\Z}{\mathbb{Z}}

\newcommand{\BX}{\mathbb{X}}
\newcommand{\BY}{\mathbb{Y}}
\newcommand{\id}{\text{id}}

\newcommand{\mMod}{\operatorname{-mod}\nolimits}

\newcommand{\C}{\mathcal{C}}
\newcommand{\CC}{\mathcal{C}}
\newcommand{\D}{\mathcal{D}}

\newcommand{\gen}[1]{\langle#1\rangle}

\newcommand{\onto}{\twoheadrightarrow}
\newcommand{\into}{\hookrightarrow}
\newcommand{\abs}[1]{\left|#1\right|}

\newcommand{\rsa}{\rightsquigarrow}
\newcommand{\arr}[1]{\stackrel{#1}{\to}}

\newcommand{\Db}{\operatorname{D^b}\nolimits}

\newcommand{\Chb}{\operatorname{Ch^b}\nolimits}

\newcommand{\rad}{\operatorname{rad}\nolimits}

\newcommand{\rsar}[1]{\stackrel{#1}{\rightsquigarrow}}

\newcommand{\End}{\operatorname{End}\nolimits}
\newcommand{\op}{{\operatorname{op}\nolimits}}

\newcommand{\ii}{\textbf{i}}
\newcommand{\jj}{\textbf{j}}

\newcommand{\DPic}{\operatorname{DPic}\nolimits}
\newcommand{\Tens}{\operatorname{Tens}\nolimits}

\newcommand{\Ch}{\Chb}
\newcommand{\perf}{\operatorname{per}\nolimits}

\newcommand{\grmodgr}{\operatorname{-grmod-}\nolimits}
\newcommand{\grmod}{\operatorname{-grmod}\nolimits}
\newcommand{\lin}{\operatorname{lin}\nolimits}
\newcommand{\proj}{\operatorname{-proj}\nolimits}

\newcommand{\grprojj}{\operatorname{-grproj-}\nolimits}
\newcommand{\dgmod}{\operatorname{-dgmod}\nolimits}
\newcommand{\dgcat}{\operatorname{dgCat}\nolimits}
\newcommand{\conv}{\operatorname{conv}\nolimits}
\newcommand{\bimcat}{A\otimes_kA^\op}
\newcommand{\grsh}[1]{\{#1\}}

\newcommand{\PreDb}{\operatorname{Pre-D^b}\nolimits}
\newcommand{\PreTr}{\operatorname{Pre-Tr}\nolimits}
\newcommand{\Tr}{\operatorname{Tr}\nolimits}
\newcommand{\Ho}{\operatorname{H_0}\nolimits}
\newcommand{\Sq}{\textbf{Sq}}
\newcommand{\Ar}{\textbf{Ar}}

\newcommand{\arbu}{\ar@{->}[ru]}
\newcommand{\argu}{\ar@{.>}[ru]}
\newcommand{\arru}{\ar@{-->}[ru]}
\newcommand{\arbr}{\ar@/^1pc/@{->}[r]}
\newcommand{\argr}{\ar@/^1pc/@{.>}[r]}
\newcommand{\arrr}{\ar@/^1pc/@{-->}[r]}
\newcommand{\arbrb}{\ar@/_1pc/@{->}[r]}
\newcommand{\argrb}{\ar@/_1pc/@{.>}[r]}
\newcommand{\arrrb}{\ar@/_1pc/@{-->}[r]}
\newcommand{\arbd}{\ar@{->}[rd]}
\newcommand{\argd}{\ar@{.>}[rd]}
\newcommand{\arrd}{\ar@{-->}[rd]}

\newcommand{\aarbu}{\ar@{->}[rru]}
\newcommand{\aargu}{\ar@{.>}[rru]}
\newcommand{\aarru}{\ar@{-->}[rru]}
\newcommand{\aarbr}{\ar@/^1pc/@{->}[rr]}
\newcommand{\aargr}{\ar@/^1pc/@{.>}[rr]}
\newcommand{\aarrr}{\ar@/^1pc/@{-->}[rr]}
\newcommand{\aarbrb}{\ar@/_1pc/@{->}[rr]}
\newcommand{\aargrb}{\ar@/_1pc/@{.>}[rr]}
\newcommand{\aarrrb}{\ar@/_1pc/@{-->}[rr]}
\newcommand{\aarbd}{\ar@{->}[rrd]}
\newcommand{\aargd}{\ar@{.>}[rrd]}
\newcommand{\aarrd}{\ar@{-->}[rrd]}

\newcommand{\brou}{\ar@{.>}[ru]}
\newcommand{\bryu}{\ar@{~>}[ru]}
\newcommand{\brvu}{\ar@{-->}[ru]}
\newcommand{\bror}{\ar@/^4pc/@{.>}[rrrr]}
\newcommand{\bryr}{\ar@/^4pc/@{~>}[rrrr]}
\newcommand{\brvr}{\ar@/^4pc/@{-->}[rrrr]}
\newcommand{\brorb}{\ar@/_4pc/@{.>}[rrrr]}
\newcommand{\bryrb}{\ar@/_4pc/@{~>}[rrrr]}
\newcommand{\brvrb}{\ar@/_4pc/@{-->}[rrrr]}
\newcommand{\brod}{\ar@{.>}[rd]}
\newcommand{\bryd}{\ar@{~>}[rd]}
\newcommand{\brvd}{\ar@{-->}[rd]}

\newtheorem{thm}{Theorem}[subsection]
\newtheorem{lem}[thm]{Lemma}
\newtheorem{prop}[thm]{Proposition}

\newtheorem{cor}[thm]{Corollary}

{\theoremstyle{definition} \newtheorem{defn}[thm]{Definition}
}
\newenvironment{pf}{Proof: \em}{ \hfill $\Box$}

\begin{document}

\maketitle

\begin{abstract}
If we have a braid group acting on a derived category by spherical twists, how does a lift of the longest element of the symmetric group act?  We give an answer to this question, using periodic twists, for the derived category of modules over a symmetric algebra.  The question has already been answered by Rouquier and Zimmermann in a special case.  We prove a lifting theorem for periodic twists, which allows us to apply their answer to the general case.

Along the way we study tensor products in derived categories of bimodules.  We also use the lifting theorem to give new proofs of two known results: the existence of braid relations and, using the theory of almost Koszul duality due to Brenner, Butler, and King, the result of Rouquier and Zimmermann mentioned above.



\end{abstract}

\tableofcontents

\setlength{\parindent}{0pt}
\setlength{\parskip}{1em plus 0.5ex minus 0.2ex}

\section{Introduction}

Actions of braid groups on derived categories were constructed by both Rouquier and Zimmermann \cite{rz} and Seidel and Thomas \cite{st} (see also \cite{hk}, \ldots).  In both cases we get a group morphism from the braid group $B_{n+1}$ on $n+1$ strands to the group of autoequivalences of a derived category $\CC$:
$$\varphi_n:B_{n+1}\to \Aut(\CC)$$
which sends the standard generators of $B_{n+1}$ to spherical twists.
In Rouquier and Zimmermann's setting, $\CC=\Db(\Gamma_n)$, where $\Gamma_n$ is the Brauer tree algebra of a line without multiplicity: this is the crucial example of an algebra whose derived category admits an action of $B_{n+1}$.

The braid group $B_{n+1}$ maps in an obvious way onto the symmetric group $S_{n+1}$ on $n+1$ letters.  $S_{n+1}$ has a unique longest element which we denote $w_0^{(n+1)}$: this is the element of $S_{n+1}$ which swaps $i$ with $n+2-i$.  We can lift this element to the braid group $B_{n+1}$.  Rouquier and Zimmermann showed that the action of the positively lifted longest element on $\Db(\Gamma_n)$ via $\varphi_n$ has a particularly nice description: it just acts by shifting complexes $n$ places to the left and twisting by an algebra automorphism.

In \cite{gra1}, a new construction of autoequivalences of derived categories of symmetric algebras, called periodic twists, was introduced, generalizing the construction of spherical twists.  In Section 6 of that article it was noted that, given an action of $B_3$ on the derived category $\Db(\Gamma_3)$, the lift of $w^0_3$ to $B_3$ acts by a periodic twist.  It is natural to ask whether this phenomenon holds in general and this question is the motivation for the current article.

Periodic twists are defined using a cone construction so it is not surprising that studying their compositions leads naturally to studying tensor products of distinguished triangles.  In Section \ref{sec-tens} we consider tensor products of distinguished triangules in derived categories of bimodules in order to set up some machinery which will be useful later.  We investigate May's axioms \cite{m} using ideas of Keller and Neeman \cite{kn} and the tools of enhanced triangulated categories due to Bondal and Kapranov \cite{bk}.

In Section \ref{sec-lift} we will prove a lifting theorem, which allows us to lift relations between periodic twists from endomorphism algebras to our original algebra: loosely, this says that relations between periodic twists that hold in the derived category of the endomorphism algebra of some projective module also hold in the derived category of the original algebra.  Combined with the Rouquier-Zimmermann description of the action of $w_0^{(n+1)}$, this will allow us to prove that positive lifts of longest elements act by periodic twists in general.  This is explained in Section \ref{sec-long}, where we also give a new, explicit, proof of Rouquier-Zimmermann's result which uses the theory of almost Koszul duality due to Brenner, Butler, and King \cite{bbk}.

An early statement of the lifting theorem was presented at the 44th Symposium on Ring Theory and Representation Theory at Okayama University in September 2011 and was included in an article submitted to the proceedings of this symposium \cite{gra2}.

\textbf{Acknowledgements:} This work was supported first by the Japan Society for the Promotion of Science and then by the Engineering and Physical Sciences Research Council [grant number EP/G007947/1].  I would like to thank Osamu Iyama for useful conversations on Koszul algebras and Robert Marsh for comments on an early version of this work.

\subsection{Conventions and notation}

All algebras are finite dimensional $k$-algebras, where $k$ is an algebraically closed field of arbitrary characteristic.  Moreover, for simplicity, all algebras are basic, i.e., their simple modules are $1$-dimensional.  By Morita theory, this is no restriction up to categorical equivalence.  Modules are finitely generated left modules unless we say otherwise. 

For an algebra $A$ over a field $k$, or over a semisimple base ring $S$, and an $A$-module $M$, we write $M^*$ to denote the dual $A$-module $\Hom_k(M,k)$ or $\Hom_S(M,S)$, as appropriate.  This duality turns left modules into right modules, and vice versa.

In Sections \ref{sec-lift} and \ref{sec-long} all algebras denoted $A$ will be symmetric, i.e, we have an isomorphism $A\cong A^*=\Hom_k(A,k)$ of $A\da A$-bimodules.  An $A$-module denoted $P$ or $P_i$ will always be projective, but will not in general be indecomposable.

We write the composition of morphisms in a category as follows: $X\arr{f}Y\arr{g}Z$ is written $g\circ f$ (and not $gf$, due to our due to our tensor product conventions, described below).  However, for the composition of arrows $1\arr{\alpha}2\arr{\beta}3$ in a quiver we write $\alpha\beta$.  

For algebras $A$ and $B$, we will say $A\da B$-bimodule when we mean $A\otimes_kB^\op$-module.  Let $\Db(A)$ and $\Db(A\da B)$ denote the derived categories of $A$-modules and $A\da B$-bimodules, respectively.  Recall that the perfect category $\per(A)$ is the full subcategory of the bounded derived category generated by the compact objects.  This means that $\per(A)$ is the full subcategory of $\Db(A)$ which consists of objects quasi-isomorphic to bounded complexes of projective $A$-modules, and $\per(A\da B)$ is the full subcategory of $\Db(A\da B)$ which consists of objects quasi-isomorphic to bounded complexes of projective $A\da B$-bimodules.  Note that the property of being a projective $A\da B$-bimodule is stronger than the property of being projective as both a left $A$-module and a right $B$-module.

Derived categories are triangulated, and morphisms of triangulated categories are triangulated functors, which take distinguished triangules to distinguished triangules and commute with the shift functor.  We denote the shift functor by $\shift$ in Section \ref{sec-tens} and by $[1]$ in Sections \ref{sec-lift} and \ref{sec-long}.  Distinguished triangles are written
$$X\to Y\to Z\to \shift X$$
or simply
$$X\to Y\to Z\rsa$$
so $Z\rsa$ is a map from $Z$ to $\shift X$.

A stalk complex is a complex where all but one constituent module is zero.  We have a full and faithful embedding $$A\mMod\into\Db(A)$$ which sends a module to the stalk complex concentrated in degree $0$, and we will use this implicitly for both modules and bimodules.

We will make much use of tensor products, but it will be useful to suppress the tensor product sign at times, especially in large commutative diagrams.  This should cause no confusion for objects as we will be clear about which categories our objects live in.  So, for example, if $P\in\Db(A\da E)$ and $Y\in\Db(E\da E)$ then we will sometimes write $PY$ instead of $P\otimes_EY$ for the object of $\Db(A\da E)$.  This is reasonable in the same way that writing a composition of functions without specifying their domains and codomains is reasonable, as the objects we tensor can be seen as $1$-morphisms in a $2$-category.  
A more uncommon notational convention that we employ is to write two juxtaposed functions to denote a tensor product, not a composition.  For example, if $f\in\End_{\Db(A\da E)}(P)$ and $g\in\End_{\Db(E\da E)}(Y)$ then $fg\in\End_{\Db(A\da E)}(PY)$.  A tensor product of categories will always be over $k$, as per Definition \ref{tenscat}.

Everything will be defined using cochain complexes $(X,d)$ which are made up of modules $X^i$ and a differential $d$ which maps $X^i\to X^{i+1}$, but as we are working in an algebraic setting and so are interested in projective resolutions, we will often write $X_i:=X^{-i}$ to avoid negative numbers.

We use the Koszul sign rule for tensor products: for $f:X\to X'$ and $g:Y\to Y'$,
$$(f\otimes g)(x\otimes y)=(-1)^{iq}f(x)\otimes g(y)$$
if $x\in X^i$ and $g$ is homogeneous of degree $q$.
We will often identify the isomorphic objects $A\otimes_AX$, $X$, and $X\otimes_BB$ in $\Db(A\da B)$ and will pretend that tensor products are strictly associative.  
This causes no serious problems.

\section{Tensor products of distinguished triangles}\label{sec-tens}

May studied symmetric tensor products of maps in distinguished triangles and gave some axioms that well-behaved triangulated categories should satisfy \cite{m}.  He showed that these axioms are satisfied in certain nice cases.  Keller and Neeman gave a nice interpretation of some of May's work based on the derived category of a commutative square and results of Happel on Dynkin quivers \cite{ha}, and noted that the symmetry of the tensor product was not of crucial importance \cite[Remark 3.11]{kn}.

We analyse the derived category of a commutative square in the spirit of Keller and Neeman, but our analysis is self-contained and does not use Happel's results.  We would like a version of \cite[Lemma 3.8]{kn} for triangulated categories.  To obtain this, we use the theory of DG-enhancements \cite{bk}.

In this section we will consider braid diagrams and braid axioms.  There is no direct connection with the braid groups mentioned in the introduction and considered in later sections.

\subsection{Modules over a commutative square}

Keller and Neeman realized that to study tensor products of distinguished triangles we can work inside the derived category of the modules over the $k$-linear category of a commutative square.  Instead of using Happel's description of this category \cite{ha} via the derived equivalence between a commutative square and an oriented graph of type $D_4$, we work directly with the commutative square.  This approach has the advantage of being self-contained, and moreover we believe it is simpler.  We will also be able to exhibit a direct connection between the derived category of the commutative square and May's braid axioms \cite{m}.

Let $\Ar$ be the path algebra of the quiver
$$Q_\Ar=1\arr{\alpha}2$$
and, for $i\in\{1,2\}$, let $P_i=\Ar e_i$ be the $i$th projective $\Ar$-module and $S_i\cong P_i/\rad P_i$ be the corresponding simple module.  Then, up to isomorphism, there are three indecomposable objects in $\Ar\mMod$: $P_1\cong S_1$, $P_2$, and $S_2$.  To save space and to produce diagrams which are more pleasing to the eye, we will write the modules $P_1$, $P_2$, and $S_2$ as simply $1$, $2$, and $3$.  Note that $P_1$ maps into $P_2$ by right multiplication by $\alpha$.  We will denote this map simply by $\alpha:P_1\into P_2$.

Our three $\Ar$-modules fit in a short exact sequence
$$0\to 1\stackrel{\alpha}{\into} 2\stackrel{\beta}{\onto} 3\to0$$
which is, up to isomorphism, the only non-split indecomposable short exact sequence in $\Ar\mMod$.  

Let $\Sq$ be the algebra $\Ar\otimes_k\Ar$.  We can write $\Sq$ as a quiver with relations in the following manner: $\Sq\cong Q_\Sq/I_\Sq$ where $Q_\Sq$ is the quiver
$$\xymatrix{
11\ar[r]^{a}\ar[d]^{c} &21\ar[d]^b \\
12\ar[r]^d &22
}$$
and $I_\Sq$ is the ideal generated by $ab-cd$.  Under the isomorphism, $a$, $b$, $c$, and $d$ correspond to $\alpha\otimes_k\id_1$, $\id_2\otimes_k\alpha$, $\id_1\otimes_k\alpha$, and $\alpha\otimes_k\id_2$ which we abbreviate as $\alpha1$, $2\alpha$, $1\alpha$, and $\alpha2$, respectively.

We can view our algebra $\Sq$ as a preadditive category with one object.  Then the $k$-linear category $\Box$ studied by Keller and Neeman is the idempotent completion of our category $\Sq$.  Hence we have equivalences
$$\Sq\mMod\cong\Box\mMod$$
of module categories and
$$\Db(\Sq)\cong\Db(\Box)$$
of derived categories.  From now on we will not use the category $\Box$ and will only consider $\Sq$.

Recall that there is a simple (naive) definition of the tensor product of two $k$-linear categories:
\begin{defn}\label{tenscat}
Given $k$-linear categories   $\C$ and $\D$, the category $\C\otimes\D$ is defined as follows: the objects of $\C\otimes\D$ are the ordered pairs $(c,d)$ for $c\in\C$ and $d\in\D$, and the formula $$\Hom_{\C\otimes\D}((c,d),(c',d'))=\Hom_\C(c,c')\otimes_k\Hom_\D(d,d')$$
gives the hom-spaces.
\end{defn}
  Then given two algebras $A$ and $B$ we have a functor
$$i_{A,B}:A\mMod\otimes B\mMod\into A\otimes_kB\mMod$$
defined in the obvious way: the tensor product of $M\in A\mMod$ and $N\in B\mMod$ is sent to $M\otimes_kN\in M\otimes_kN\mMod$.
\begin{lem}\label{fftens}
 For $M_i\in A\mMod$ and $N_i\in B\mMod$, $i=1,2$, there is a natural isomorphism of vector spaces
$$\Hom_{A\otimes_k B}(M_1\otimes_k N_1,M_2\otimes_k N_2)\cong \Hom_A(M_1,M_2)\otimes_k\Hom_B(N_1,N_2)\fs$$

\begin{pf}
 Recall the definition: for an algebra $\Lambda$ and $\Lambda$-modules $V$ and $W$,
$$\Hom_\Lambda(V,W)=\left\{f\in\Hom_k(V,W)\:|\:f(\lambda v)=\lambda f(v) \fa \lambda\in\Lambda, v\in V\right\}\subset\Hom_k(V,W)\fs$$
Use the isomorphism $\Hom_k(V,W)\cong V^*\otimes_k W$ which sends $f$ to $\sum_i v_i^*\otimes f(v_i)$ for some basis $\{v_i\}$ of $V$.  Then $\Hom_\Lambda(V,W)$ corresponds to the subspace 
$$\left\{\sum g_i\otimes w_i\in V^*\otimes_k W\:|\:  \sum g_i\lambda\otimes w_i=\sum g_i\otimes \lambda w_i \fa \lambda\in\Lambda\right\}\subset V^*\otimes_k W\fs$$
Using the braiding in the symmetric monoidal category of vector spaces gives an isomorphism
$$(M_1\otimes_k N_1)^*\otimes_k (M_2\otimes_k N_2)\cong N_1^*\otimes_k M_1^*\otimes_k M_2\otimes_k N_2 \cong (M_1^*\otimes_k M_2)\otimes_k (N_1^*\otimes_k N_2 )$$
and we can check that the subspace corresponding to $\Hom_{A\otimes_k B}(M_1\otimes_k N_1,M_2\otimes_k N_2)$ and the subspace corresponding to $\Hom_A(M_1,M_2)\otimes_k\Hom_B(N_1,N_2)$ agree.  Naturality is clear.
\end{pf}
\end{lem}
So the functor $i_{A,B}$ is fully faithful.  Note that in general $A\mMod\otimes B\mMod$ is not an abelian category.  It will be abelian, and the functor will in fact be an equivalence, if and only if at least one of $A$ and $B$ is semisimple.

Using the above functor we can immediately write down nine non-isomorphic indecomposable $\Sq$-modules: this just comes down to tensoring one of the $\Ar$-modules labelled $1$, $2$, or $3$ with another.  We see that they fit in the $3\times 3$ commutative diagram of short exact sequences
$$\xymatrix{
 & 0\ar[d] & 0\ar[d] & 0\ar[d] & \\
0\ar[r] & 11\ar[r]^{\alpha1}\ar[d]^{1\alpha} & 21\ar[r]^{\beta1}\ar[d]^{2\alpha} & 31\ar[r]\ar[d]^{3\alpha} & 0 \\
0\ar[r] & 12\ar[r]^{\alpha2}\ar[d]^{1\beta} & 22\ar[r]^{\beta2}\ar[d]^{2\beta} & 32\ar[r]\ar[d]^{3\beta} &0 \\
0\ar[r] & 13\ar[r]^{\alpha3}\ar[d] & 23\ar[r]^{\beta3}\ar[d] & 33\ar[r]\ar[d] & 0\\
 &0 &0 &0 &
}$$
of $\Sq$-modules, where we have suppressed the tensor product sign.

To list every object in $\Sq\mMod$, we can take the abelian hull of $\Ar\mMod\otimes\Ar\mMod$ in $\Sq\mMod$: we start with our list $11,21,\ldots,33$ and add the kernel and cokernel of every map between these modules.  It is easy to see that the only modules that are missing from our list are the the kernel of the epimorphism $\beta\otimes\beta:22\onto 33$ and the cokernel of the monomorphism $\alpha\otimes\alpha:11\into 22$, which we denote $K$ and $C$, respectively.  So $\Sq$ has finite representation type, and in fact its eleven modules fit into a large commutative diagram
$$
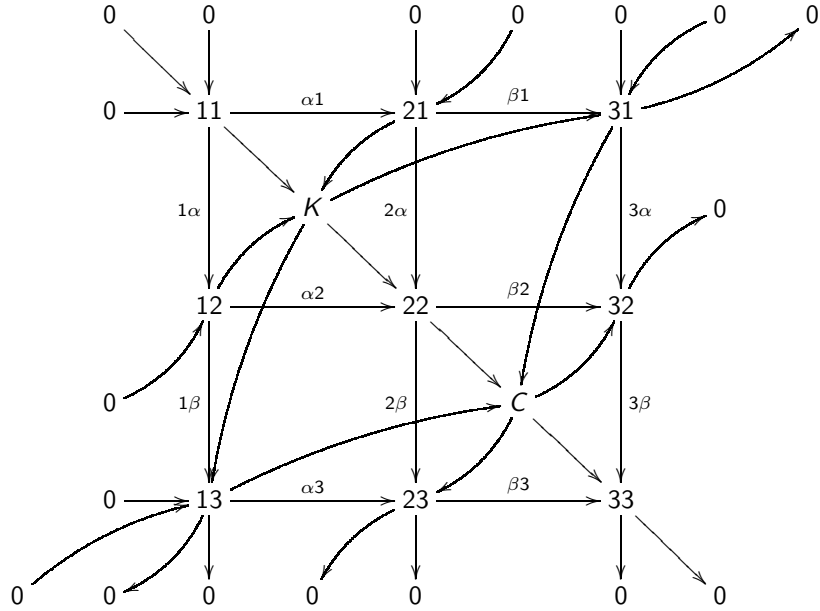
\xymatrix{
  & 0\ar[rd] &0\ar[d] & & 0\ar[d] &0\ar@/^/[ld] &0\ar[d] & 0\ar@/_/[ld] &0\\
&0\ar[r]&11\ar[rr]^{\alpha 1}\ar[dd]_{1\alpha}\ar[rd] & & 21\ar[rr]^{\beta  1}\ar[dd]_{2 \alpha}\ar@/_/[ld] & & 31\ar@/_/[dddl]   \ar[dd]^{3 \alpha} \ar@/_/[rru] & & \\
&& & K\ar@/^/[rrru]^{}\ar@/_/[dddl]^{}\ar[rd] &&&&0&\\
&&12\ar[rr]^{\alpha  2}\ar[dd]_{1 \beta}\ar@/^/[ur] & &
22\ar[rd]^{}\ar[dd]_{2 \beta}\ar[rr]^{\beta  2} &  & 32\ar[dd]^{3 \beta}\ar@/^/[ur] & & \\
&0\ar@/_/[ur]  &&& &C\ar@/^/[ld] \ar@/_/[ur]^{ }\ar[rd]  &&&\\
&0\ar[r]&13\ar[rr]^{\alpha  3}\ar@/^/[rrru]^{ }\ar@/^/[ld]\ar[d] &&
23\ar[rr]^{\beta 3}\ar[d]\ar@/_/[ld] & & 33\ar[d]\ar[rd]&&\\
0\ar@/^/[rru] &0&0&0&0&&0&0&
}$$
consisting of many short exact sequences.  But we will see that the situation is even clearer on the derived level.  Note that, as well as the short exact sequences that can be read off from the diagram, we have a short exact sequence
$$0\to 11\into 12\oplus 21 \onto K\to 0$$
which we call the \emph{Mayer-Vietoris} short exact sequence.

The maps between modules, if not all the short exact sequences, are more clearly read off from the Auslander-Reiten quiver of $\Sq\mMod$:
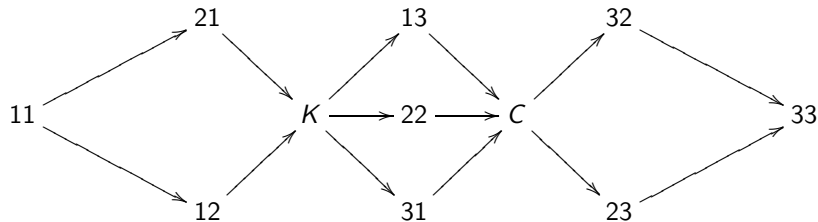
$$ \xymatrix{
& & 21\ar[rd] & & 13\ar[rd] & & 32\ar[rrd] & &  & \\
11\ar[rru]\ar[rrd] & && K\ar[r]\ar[ur]\ar[dr] & 22\ar[r] & C\ar[ru]\ar[rd] &  &  & 33 \\
& & 12\ar[ru] & & 31\ar[ru] & & 23\ar[rru] & &  &
} $$
Twisting this quiver and curving its arrows gives us the diagram
$$ 
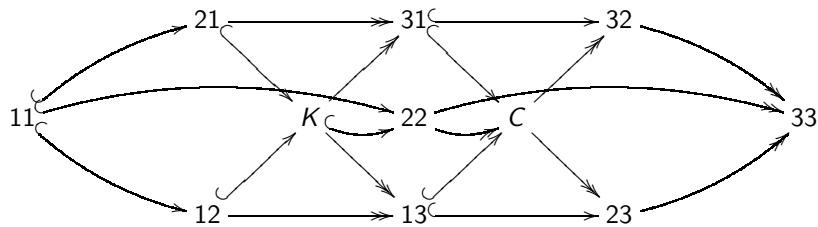
\xymatrix{
& & 21\ar@{^{(}->}[rd]\ar@{->>}[rr] & & 31\ar@{^{(}->}[rd]\ar@{^{(}->}
[rr] & & 32\ar@{->>}@/^/[rrd] & &  & \\
11\ar@{^{(}->}@/^/[rru]\ar@{^{(}->}@/_/[rrd]\ar@{^{(}->}@/^1pc/[rrrr] & && K\ar@{^{(}->}@/_/[r]\ar@{->>}[ur]\ar@{->>}[dr] & 22\ar@{->>}@/_/[r]\ar@{->>}@/^1pc/[rrrr] & C\ar@{->>}[ru]\ar@{->>}[rd] &  &  & 33 \\
& & 12\ar@{^{(}->}[ru]\ar@{->>}[rr] & & 13\ar@{^{(}->}[ru]\ar@{^{(}->}[rr] & & 23\ar@{->>}@/_/[rru] & &  &
} $$
which better illustrates the short exact sequences.  The reader should imagine the arrows of the short exact sequences
$$0\to11\into22\onto C\to0$$
and
$$0\to K\into22\onto 33\to0$$
as being threaded in and out of the diagram, going into and out of the piece of paper or computer screen.

\subsection{The derived category of a commutative square}

As $\Ar$ is a path algebra with no relations, it is hereditary and so has global dimension $1$.  As $\Ar$ is hereditary, the three objects $1,2,3$ are, up to isomorphism and shift, the only indecomposable objects in $\Db(\Ar)$, and
$$1\arr{\alpha} 2\arr{\beta} 3\rsar{\gamma}$$
is, up to rotation and isomorphism, the only indecomposable distinguished triangle.  Note that we have labelled the map $3\to \shift 1$ as $\gamma$.

Using the identification of $\Ar$ with upper triangular $2\times 2$ matrices with entries in $k$, we can identify $\Sq$ with upper triangular matrices with entries in $\Ar$, and so conclude \cite{erz} that $\Sq$ has global dimension $2$.  This is also easy to see by directly calculating the projective resolutions of all simple modules.  

We can check that the cones of the maps $\alpha\otimes\beta:12\to 23$ and $\beta\otimes\alpha:21\to 32$ are quasi-isomorphic and are not concentrated in degree zero.  Let $G$ denote a representative of the isomorphism class of objects containing the two cones above in the derived category $\Db(\Sq)$.  As $\Sq$ has global dimension $2$, all isomorphism classes of indecomposable objects in $\Db(\Sq)$ have a representative given by a chain complex of projective $\Sq$-modules with at most three nonzero terms, so we can quicky check that, up to isomorphism and shift, this is the only object of $\Db(\Sq)$ which is not concentrated in a single degree.  So we have found all twelve isoclasses of indecomposable objects, up to shift, in $\Db(\Sq)$.

\begin{lem}\label{tenstri}
Suppose $A$ and $B$ are $k$-algebras and we have an $A$-module $M$ and a distinguished triangle
$$X\arr{f} Y\arr{g} Z\arr{h}\shift X$$
in $\Db(B)$.  Then the triangle defined by the bottom row of the commutative diagram
$$\xymatrix{
M\otimes_k X\ar[r]^{M\otimes f}\ar@{=}[d] & M\otimes_k Y\ar[r]^{M\otimes g}\ar@{=}[d] & M\otimes_k Z\ar[r]^{M\otimes h}\ar@{=}[d] & M\otimes_k (\shift X)\ar[d]^{\sim}\\
M\otimes_k X\ar[r]^{M\otimes f} & M\otimes_k Y\ar[r]^{M\otimes g} & M\otimes_k Z\ar[r]^{} & \shift(M\otimes_k X)
}$$
is distinguished in $\Db(A\otimes_k B)$.

\begin{pf}
Consider the canonical distinguished triangle
$$M\otimes_k X\arr{M\otimes f} M\otimes_k Y\arr{} \cone(M\otimes f) \arr{} \shift(M\otimes_k X)$$
in the homotopy category.  By considering the first map in this triangle, it is clearly isomorphic to the triangle displayed in the lemma.  Now localize to pass to the derived category.
\end{pf}
\end{lem}

As in \cite[Remark 4.2]{m}, we will implicitly use fixed isomorphisms
$$(\shift i)\otimes j\cong\shift(i\otimes j)\cong i\otimes (\shift j)$$
for $i,j\in\{1,2,3\}$.
These come from the isomorphisms
$$(\shift X)\otimes Y\cong\shift(X\otimes Y)\cong X\otimes (\shift Y)$$
$$(-1)^{\deg(y)}x\otimes y\mapsto x\otimes y\mapsto (-1)^{\deg(x)}x\otimes y$$
of differential graded $k$-modules.
The notation $\shift{ij}$ will mean $\shift{(i\otimes j)}$.

We will sometimes be sloppy and write maps such as
$$31\arr{\gamma1}\shift(11)$$
when we really mean the composition
$$31\arr{\gamma1}(\shift 1)1\arr{\sim}\shift(11)$$
for our fixed isomorphism $(\shift 1)1\arr{\sim}\shift(11)$.

Now we must be very careful with our signs.  
We consider the diagram
$$\xymatrix @!=5pt {
11\ar[rr]^{\alpha   1}\ar[dd]^{1  \alpha} & & 21\ar[rr]^{\beta   1}\ar[dd]^{2  \alpha} & & 31\ar[rr]^{\gamma   1}\ar[dd]^{3  \alpha} & & \shift 11\ar[dd]^{\shift 1\alpha}\\
 &&&&\\
12\ar[rr]^{\alpha   2}\ar[dd]^{1  \beta} & &
22\ar[dd]^{2  \beta}\ar[rr]^{\beta   2} &  &  32\ar[rr]^{\gamma   2}\ar[dd]^{3  \beta} & & \shift 12\ar[dd]^{\shift 1  \beta}\\
\\
13\ar[rr]^{\alpha   3}\ar[dd]^{1  \gamma}
&&
23\ar[rr]^{\beta  3}\ar[dd]^{2  \gamma} & & 33\ar[rr]^{\gamma   3}\ar[dd]^{3  \gamma} & & \shift 13\ar[dd]^{-\shift 1  \gamma}\\
&&&&&(-) \\
\shift 11\ar[rr]^{\shift \alpha   1} & & \shift 21\ar[rr]^{\shift \beta   1} & & \shift 31\ar[rr]^{-\shift \gamma   1} & & \shift^2 11
}$$
in $\Db(\Sq)$.  By Lemma \ref{tenstri} and the axioms for a triangulated category, each row and column is a distinguished triangle, and we can check directly that each square commutes except the bottom right square (marked $(-)$) which anticommutes.  The diagram should be compared to Lemma 2.6 of \cite{m}. 
This can also be viewed as a braid of morphisms of triangles, as follows:
$$
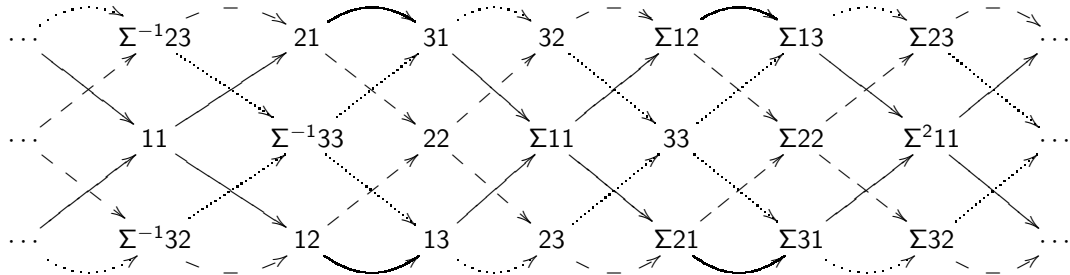
\xymatrix{
\cdots\argr\arbd &\shift^{-1} 23\arrr\argd &21\arbr\arrd &31\argr\arbd &32\arrr\argd &\shift 12\arbr\arrd &\shift 13\argr\arbd &\shift 23\arrr\argd &\cdots \\
\cdots\arru\arrd &11\arbu\arbd &\shift^{-1}33\argu\argd &22\arru\arrd &\shift 11\arbu\arbd &33\argu\argd &\shift 22\arru\arrd &\shift^2 11\arbu\arbd &\cdots \\
\cdots\arbu\argrb &\shift^{-1}32\argu\arrrb &12\arru\arbrb &13\arbu\argrb &23\argu\arrrb &\shift 21\arru\arbrb &\shift 31\arbu\argrb &\shift 32\argu\arrrb &\cdots
}$$
where all quadrilaterals
commute except the square
$$
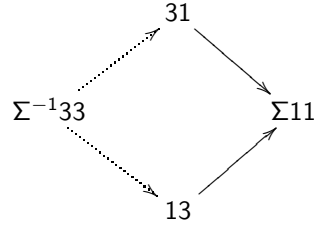
\xymatrix{
&31\arbd &\\
\shift^{-1}33\argu\argd &&\shift 11  \\
&13\arbu &
}$$
and its shifts, which anticommute.
These ``braid diagrams'' appear in, for example, \cite{ive,m}.  
We have used different types of arrows to help the reader see the different distinguished triangles.

We obtain a more complicated diagram
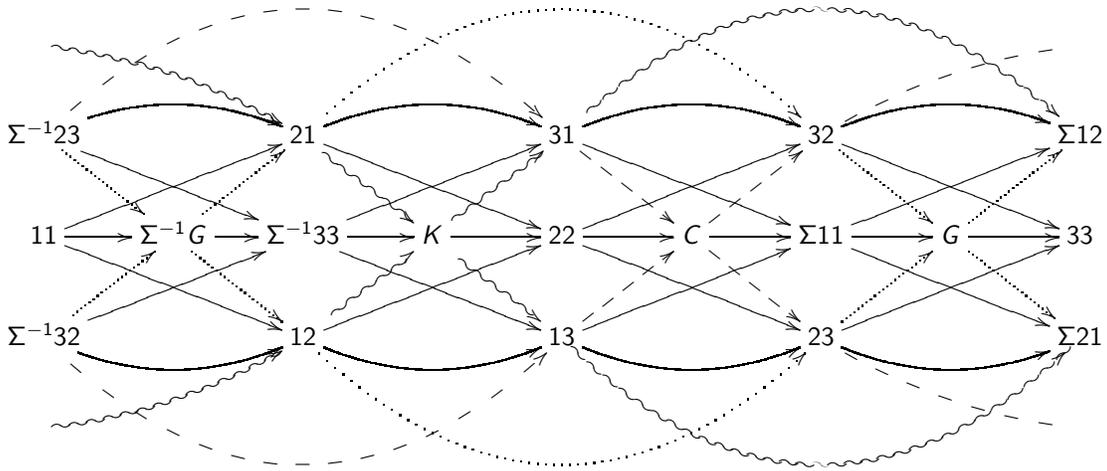
\vspace{1em}
$$\xymatrix @!C=25pt {
\ar@/^/@{~>}[rrd]&&&&&&&&\\
\shift^{-1}23\ar@/^1pc/[rr]\ar[rrd]\brod\brvr &&21\ar@/^1pc/[rr]\ar[rrd]\bror\bryd &&31\ar@/^1pc/[rr]\ar[rrd]\bryr\brvd &&32\ar@/^1pc/[rr]\ar[rrd]\brod\ar@/^/@{--}[urr]
 &&\shift12
 \\
11\ar@{->}[rru]\ar[rrd]\ar[r] &\shift^{-1}G\ar[r]\brou\brod &\shift^{-1}33\ar[r]\ar[rru]\ar[rrd] &K\ar[r]\bryu\bryd &22\ar[r]\ar[rru]\ar[rrd] &C\ar[r]\brvu\brvd &\shift11\ar[r]\ar@{->}[rru]\ar[rrd] &G\ar[r]\brou\brod &33
 \\
\shift^{-1}32\ar[rru]\ar@/_1pc/[rr]\brou\brvrb &&12\ar[rru]\ar@/_1pc/[rr]\bryu\brorb &&13\ar@{->}[rru]\ar@/_1pc/[rr]\brvu\bryrb &&23\ar[rru]\ar@/_1pc/[rr]\brou\ar@/_/@{--}[rrd]
 &&\shift21
\\
{}\ar@/_/@{~>}[rru] &&&&&&&&\\
&&&&&&&&
}$$
if we include our extra objects $K$, $C$, and $G$ and six new distinguished triangles, such as
$$21\to K\to 13\rsa$$
and
$$12\to K\to 31\rsa$$
for $K$.
Again, we have used different types of arrows to help the reader see the distinguished triangles. 

It is important to note that this diagram does not commute, even up to sign.  However, every subdiagram containing at most one of the objects $K$, $C$, and $G$ will commute up to sign. 

We note that the distinguished triangles
$$11\arr{} 22\arr{} C\rsar{}\cma$$
$$K\arr{} 22\arr{} 33\rsar{}\cma$$
and
$$\shift11\arr{}G\arr{} 33\rsar{}$$
are not visible in this diagram.  Similarly to the case of the module category, one can imagine these distinguished triangles threading in and out of the diagram, lying in a plane which is perpendicular to the plane in which the rest of the diagram is embedded.

The above diagram is just a combination of May's braid axioms (TC3) for the object $K$ and (TC3') for the object $C$ \cite[Section 4]{m}, together with the analogous braid for the object $G$.

Removing all but the irreducible morphisms from the above diagram gives the following graph:
$$ 
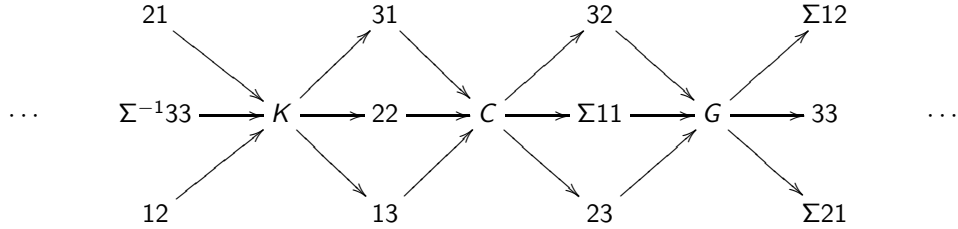
\xymatrix{
 & 21\ar[rd] & & 31\ar[rd] & & 32\ar[rd] & & \shift12 & \\
\cdots & \shift^{-1}33\ar[r] & K\ar[r]\ar[ur]\ar[dr] & 22\ar[r] & C\ar[ru]\ar[r]\ar[rd] & \shift11\ar[r] & G\ar[ru]\ar[r]\ar[rd] & 33 & \cdots\\
& 12\ar[ru] & & 13\ar[ru] & & 23\ar[ru] & & \shift21 &
} $$
which one can see is a twisted version of the Auslander-Reiten quiver
$$ 
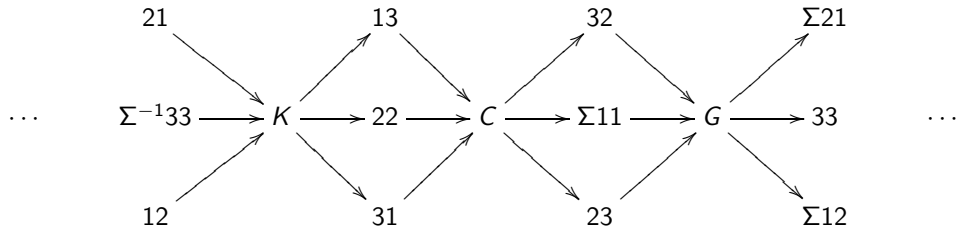
\xymatrix{
 & 21\ar[rd] & & 13\ar[rd] & & 32\ar[rd] & & \shift21 & \\
\cdots & \shift^{-1}33\ar[r] & K\ar[r]\ar[ur]\ar[dr] & 22\ar[r] & C\ar[ru]\ar[r]\ar[rd] & \shift11\ar[r] & G\ar[ru]\ar[r]\ar[rd] & 33 & \cdots\\
& 12\ar[ru] & & 31\ar[ru] & & 23\ar[ru] & & \shift12 &
} $$
of $\Db(\Sq)$.

We note that there are other distinguished triangles in $\Db(\Sq)$, the most important of which for us is the Mayer-Vietoris distinguished triangle
$$11\to 12\oplus21 \to K\rsar{}$$
coming from the associated short exact sequence in $\Sq\mMod$.
We refer the reader to \cite{kn} for more information on the distinguished triangles we have not described here.

\subsection{Getting a morphism of triangulated categories}

We want to show that the results we obtained in the last subsection hold for arbitrary tensor products in derived categories.  We will use the idea of Keller and Neeman \cite{kn}: if we have a map of triangulated categories from $\Db(\Sq)$ to another derived category then because this map sends distinguished triangules to distinguished triangules we will know that the images of the objects in $\Sq\mMod$ fit into some nice distinguished triangules.  Our objects in the derived category may not be concentrated in degree $0$ so we approach the problem by working with enhancements of triangulated categories, following Bondal and Kapranov \cite{bk}.  This strategy is more technical than, though very similar in spirit to, Example 3.2 of \cite{kn}.

Our reference throughout is \cite{bk}.

Let $\CC$ be a DG-category and let $\CC^\oplus$ denote the DG-category obtained from $\CC$ by adjoining finite formal direct sums of objects.
\begin{defn}[Bondal-Kapranov] A \emph{twisted complex over $\CC$} is a set $\{(E_i)_{i\in\Z}, q_{ij}:E_i\to E_j\}$ where the $E_i$ are objects in $\CC^\oplus$, equal to $0$ for almost all $i$, and the $q_{ij}$ are morphisms in $\CC$ of degree $i-j+1$
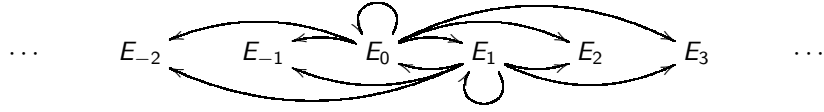
\hspace{0.5em}
$$\xymatrix{
 &&&&&&\\
{}\cdots &E_{-2} &E_{-1} &E_0\ar@/_1pc/[ll]\ar@/_0.5pc/[l]\ar@(ur,ul)[]\ar@/^0.5pc/[r]\ar@/^1pc/[rr]\ar@/^1.5pc/[rrr] &E_1\ar@/_1pc/[rr]\ar@/_0.5pc/[r]\ar@(dr,dl)[]\ar@/^0.5pc/[l]\ar@/^1pc/[ll]\ar@/^1.5pc/[lll] &E_2 &E_3 &{}\cdots\\
 &&&&&&
}$$
satisfying the Maurer-Cartan equation
 $$dq_{ij}+\sum_k q_{kj}q_{ik} =0\fs$$
Twisted complexes over $\CC$ form a DG-category, denoted $\PreTr(\CC)$, with morphism complexes made up of spaces
$$\Hom^k_{\PreTr(\CC)}(K,K')=\coprod_{i,j\in\Z}\Hom^{k+i-j}_\C(K_i,K'_j)$$
and a differential which is described in Section 1 of \cite{bk}.  Taking the category of twisted complexes is an endofunctor
$$\PreTr:\dgcat\to\dgcat$$
of the category of DG-categories.  We call it the \emph{pretriangulated completion} functor.
\end{defn}

Any DG-category has an associated homology category $H(\CC)$ with the same objects as $\CC$ but with hom spaces defined as the homology of the differential.  In particular, we can take the $0$th homology $\Ho(\CC)$ of $\CC$.  We denote $\Ho(\PreTr(\CC))$ by $\Tr(\C)$.  It is a triangulated category \cite[Section 1, Proposition 1]{bk}.

Recall that the DG-category of chain complexes of vector spaces has as morphisms $f:X\to Y$ of degree $n$ all collections of maps $\{f_i:X_i\to Y_{i-n}\}$ with no requirement that these maps commute with the differentials of $X$ and $Y$.  The differential $d$ on this DG category is such that the degree zero maps $f$ which satisfy $df=0$ are exactly the traditional maps of chain complexes, which do commute with the differentials of $X$ and $Y$.
Then, for any DG-category $\C$, a DG $\C$-module is a functor from $\C$ to the DG-category of chain complexes of vector spaces. 
We denote the category of DG $\C$-modules by $\C\dgmod$.

There is an embedding $$\iota_\C:\PreTr(\CC)\into\CC\dgmod$$ defined as follows.  Let $E=\{E_i,q_{ij}\}$ be a twisted complex over $\C$ and let $c\in\C$.  For each $i\in\Z$, write the chain complex $\Hom_{\C^\oplus}(c,E_{i})$ as
$$\cdots\arr{d_{i}} \Hom_{\C^\oplus}(c,E_{i})_{-1}\arr{d_i} \Hom_{\C^\oplus}(c,E_{i})_0\arr{d_i} \Hom_{\C^\oplus}(c,E_{i})_1\arr{d_i}\cdots$$
with differential $d_i$ of degree $1$.  
Then 
$\iota_\C(E)$ is the functor which sends $c\in\C$ to the chain complex
$$\cdots\arr{\partial_E} \bigoplus_{i+j=-1}\Hom_{\C^\oplus}(c,E_{i})_j\arr{\partial_E} \bigoplus_{i+j=0}\Hom_{\C^\oplus}(c,E_{i})_j\arr{\partial_E} \bigoplus_{i+j=1}\Hom_{\C^\oplus}(c,E_{i})_j\arr{\partial_E}\cdots$$
where $$\partial_E|_{\Hom_{\C^\oplus}(c,E_{i})_j}=d_i|_{\Hom_{\C^\oplus}(c,E_{i})_j}+\sum_{k\in\Z}\Hom_{\C^\oplus}(c,q_{ik})\fs$$
If $f\in\Hom^\ell_\C(E_i,E'_j)$ is a morphism of twisted complexes from $E=\{E_i,q_{ij}\}$ to $E'=\{E'_i,q'_{ij}\}$ then $\iota_\C(f)$ is the obvious degree $\ell$ map of chain complexes, which does not necessarily commute with the differentials $\partial_E$ and $\partial_{E'}$.

\begin{defn}[Bondal-Kapranov]
 $\C$ is called \emph{pretriangulated} if the image of every twisted complex under the above embedding is a representable functor.  If $\C$ is pretriangulated then every object in $\PreTr\C$ has a \emph{convolution}, which is defined as the associated representing object of $\C$.  A DG-functor $F:\C\to\D$ between pretriangulated categories is \emph{pre-exact} if it commutes with the operation of taking convolutions of twisted complexes, i.e., the diagram
$$\xymatrix @C=40pt {
\PreTr(\C)\ar[r]^{\PreTr(F)}\ar[d]^{\conv_{\C}} &\PreTr(\D)\ar[d]^{\conv_\D}\\
\C\ar[r]^F &\D
}$$
commutes, where $\conv$ denotes the operation of taking the convolution of a twisted complex.
\end{defn}
Note that if $F$ is a pre-exact functor then $\Ho(F)$ is a triangulated functor \cite[Section 3]{bk}.

The following definition is based on, but different to, that of \cite{bk}: we use projectives instead of injectives because we have algebraic and not geometric applications in mind.  But by the usual duality everything still works.
\begin{defn}
 Let $A$ be an algebra and let $\C$ be the DG-category obtained by treating $A$ as a preadditive category with one object.  Let $\PreDb(\C)$ be the full DG-subcategory of chain complexes of 
$\C$-modules consisting of 
complexes where each module is projective over $A$ and only finitely many modules in the complex are non-zero.
\end{defn}
$\PreDb(A)$ is a \emph{DG-enhancement} of $\Db(A)$, i.e., there is an equivalence of triangulated categories $\Ho(\PreDb(A))\arr{\sim}\Db(A)$ \cite[Section 3, Example 3]{bk}.

The following lemma states that the convolution for $\PreTr(A)$ is particularly simple: it is a generalized cone construction.  This is surely well-known to the experts.  The proof is immediate from the definitions.
\begin{lem}\label{convpredb}
If $E=\{E_i,q_{ij}\}$ is a twisted complex over $\PreDb(A)$ then its convolution is the chain complex
$$\cdots\arr{\partial} \bigoplus_{i+j=-1}E_{i,j}\arr{\partial} \bigoplus_{i+j=0}E_{i,j}\arr{\partial} \bigoplus_{i+j=1}E_{i,j}\arr{\partial}\cdots$$
in $\PreDb(A\da A)$, where $E_i$ is the complex
$$\cdots\arr{d_{i}}E_{i,-1}\arr{d_i} E_{i,0}\arr{d_i}E_{i,1}\arr{d_i}\cdots$$
and $$\partial|_{\bigoplus_{i+j=k}E_{i,j}}=d_i|_{E_{i,j}}+\sum_{k\in\Z}q_{ik}\fs$$
\end{lem}

Now we can start laying the first few bricks in the construction of our functor.

\begin{lem}\label{constr}
 Let $f:X\to Y$ be a morphism in $\Db(A\da A)$.  Then there is a pre-exact morphism of pretriangulated categories
$$F:\PreDb(\Ar)\to\PreDb(A\da A)$$
such that, on taking homotopy, the map
$$\Ho(F):\Db(\Ar)\to\Db(A\da A)$$
is a triangulated functor which, up to isomorphism, takes the morphism $\alpha:P_1\to P_2$ of stalk complexes to $f:X\to Y$.

\begin{pf}
We take projective $A\da A$-bimodule resolutions $\BX$ and $\BY$ of $X$ and $Y$ and lift $f$ to a map $\varphi:\BX\to \BY$ of chain complexes.  Any object $Z$ in $\PreDb(\Ar)$ is a complex made up of direct sums of $P_1$ and $P_2$ with differentials only consisting of multiples of $\alpha$ and identity maps.  There is a functor $G:\Ar\proj\to\PreDb(A\otimes_kA^\op)$ defined on the additive category of projective $\Ar$-modules which sends $P_1$, $P_2$, and $\alpha$ to $\BX$, $\BY$, and $\varphi$, respectively.  So $G$ induces a functor from bounded chain complexes of projective $\Ar$-modules into chain complexes in the category $\PreDb(A\otimes_kA^\op)$, and taking the total complex gives us a functor $F:\PreDb(\Ar)\to\PreDb(A\otimes_kA^\op)$.  In other words, $F$ is defined by a cone construction.
As we are working with chain complexes, this is functorial.

It is clear that $F$ is a DG-morphism and we now check that it is pre-exact, i.e., that the diagram
$$\xymatrix @C=50pt {
\PreTr(\PreDb(\Ar))\ar[r]^{\PreTr(F)}\ar[d]^{\conv_{\PreDb(\Ar)}} &\PreTr(\PreDb(\bimcat))\ar[d]^{\conv_{\PreDb(\bimcat)}}\\
\PreDb(\Ar)\ar[r]^F &\PreDb(\bimcat)
}$$
commutes.  This comes down to a routine check which we will now describe.

Let $E=(E_i,q_{ij})\in\PreTr(\PreDb(\Ar))$, and for $M\in\Ar\proj$, write $G(M)$ as
$$\cdots\to G_{-1}(M)\to G_{0}(M)\to G_{1}(M)\to\cdots\fs$$

Using Lemma \ref{convpredb}, the convolution of $E$ is 
$$\cdots\arr{\partial} \bigoplus_{i+j=-1}E_{i,j}\arr{\partial} \bigoplus_{i+j=0}E_{i,j}\arr{\partial} \bigoplus_{i+j=1}E_{i,j}\arr{\partial}\cdots$$
and then to find the image of this complex under $F$ we take the total complex of
$$\xymatrix{
 &\vdots\ar[d] &\vdots\ar[d] &\vdots\ar[d] &\\
\cdots\ar[r] &\bigoplus_{i+j=-1}G_{-1}(E_{i,j})\ar[r]\ar[d] &\bigoplus_{i+j=0}G_{-1}(E_{i,j})\ar[r]\ar[d] &\bigoplus_{i+j=1}G_{-1}(E_{i,j})\ar[r]\ar[d] & \cdots \\
\cdots\ar[r] &\bigoplus_{i+j=-1}G_{0}(E_{i,j})\ar[r]\ar[d] &\bigoplus_{i+j=0}G_{0}(E_{i,j})\ar[r]\ar[d] &\bigoplus_{i+j=1}G_{0}(E_{i,j})\ar[r]\ar[d] & \cdots \\
\cdots\ar[r] &\bigoplus_{i+j=-1}G_{1}(E_{i,j})\ar[r]\ar[d] &\bigoplus_{i+j=0}G_{1}(E_{i,j})\ar[r]\ar[d] &\bigoplus_{i+j=1}G_{1}(E_{i,j})\ar[r]\ar[d] & \cdots \\
 &\vdots &\vdots &\vdots &
}$$
which gives the complex
$$\cdots\arr{} \bigoplus_{i+j+k=-1}G_{k}(E_{i,j})\arr{} \bigoplus_{i+j+k=0}G_{k}(E_{i,j})\arr{} \bigoplus_{i+j+k=1}G_{k}(E_{i,j})\arr{}\cdots$$
in $\PreDb(\bimcat)$.

Now let's go the other way around the square.  Applying $\PreTr(F)$ to $E$ just gives the twisted complex $\{F(E_i),F(q_{ij})\}$.  Unwinding the definitions, we see that $F(E_i)$ is the complex
$$\cdots \bigoplus_{j+k=-1}G_k(E_{i,j})\to \bigoplus_{j+k=0}G_k(E_{i,j})\to \bigoplus_{j+k=1}G_k(E_{i,j})\to \cdots$$
So, using Lemma \ref{convpredb} again, the convolution of $\{F(E_i),F(q_{ij})\}$ is
$$\cdots\arr{} \bigoplus_{i+j=-1}F(E_i)_j\arr{} \bigoplus_{i+j=0}F(E_i)_j\arr{} \bigoplus_{i+j=1}F(E_i)_j\arr{}\cdots$$
which we can rewrite as
$$\cdots\arr{} \bigoplus_{i+j+k=-1}G_{k}(E_{i,j})\arr{} \bigoplus_{i+j+k=0}G_{k}(E_{i,j})\arr{} \bigoplus_{i+j+k=1}G_{k}(E_{i,j})\arr{}\cdots$$
and we see that this has the same terms as the complex above.  It is simple to check that the differentials also agree.

So $F$ is a pre-exact morphism and hence $\Ho(F)$ is a triangulated functor.  By definition it takes $\alpha$ to $\varphi$, which is isomorphic to $f$.
\end{pf}
\end{lem}

We will need to find a way to tensor together two of the maps constructed using the previous lemma.  The following result will be useful.
\begin{lem}\label{tenspre}
For two algebras $A,B$, the obvious inclusion
$$i_{A,B}:\PreDb(A)\otimes\PreDb(B)\to\PreDb(A\otimes_kB)$$
becomes a quasi-equivalence
$$\PreTr(i_{A,B}):\PreTr(\PreDb(A)\otimes\PreDb(B))\to\PreTr(\PreDb(A\otimes_kB))$$
on applying the pretriangulated completion functor, i.e., the functor
$$\Ho(\PreTr(i_{A,B})):\Tr(\PreDb(A)\otimes\PreDb(B))\to\Tr(\PreDb(A\otimes_kB))\cong\Db(A\otimes_kB)$$
is an equivalence of triangulated categories. 

\begin{pf}
If $\C$ is a pretriangulated category then $\Ho(\PreTr(\C))$ is equivalent to $\Ho(\C)$ \cite[Section 3, Proposition 1]{bk}, so the codomain of our functor
$$\Ho\left(\PreTr(\PreDb(A)\otimes_k\PreDb(B))\right)\to\Ho\left(\PreTr(\PreDb(A\otimes_kB))\right)$$
is equivalent to $\Db(A\otimes_k B)$.  The domain is the $0$th homology of a pretriangulated category, and so is triangulated.

By extending Lemma \ref{fftens} we see that our our original functor $i_{A,B}$ is fully faithful.  From the definition, $\PreTr(i_{A,B})$ and therefore $\Ho(\PreTr(i_{A,B}))$ are also both fully faithful, and so the domain of our functor is isomorphic to a full subcategory of $\Db(A\otimes_k B)$.  But it clearly contains all simple $A\otimes_k B$-modules and so, as it is closed under taking cones, must be the whole subcategory.
\end{pf}
\end{lem}

We can now construct our morphism.
\begin{prop}\label{thetrimap}
For $i\in\{1,2\}$, let $F_i:\PreTr(\Ar)\to\PreTr(\bimcat)$ be exact morphisms of pretriangulated categories.  Then there is a triangulated functor
$$F:\Db(\Sq)\to\Db(\bimcat)$$
which acts as $\Ho(F_1)\otimes\Ho(F_2)$ on the subcategory $\Db(\Ar)\otimes\Db(\Ar)$ of $\Db(\Sq)$.

\begin{pf}
Tensor together our two functors $F_i$, $i\in\{1,2\}$, to get a morphism
$$F_1\otimes_k F_2:\PreDb(\Ar)\otimes\PreDb(\Ar)\to\PreDb(\bimcat)\otimes\PreDb(\bimcat)$$
of dg-categories.  We compose this with the tensor product functor $$\PreDb(\bimcat)\otimes\PreDb(\bimcat)\to\PreDb(\bimcat)$$ to get a morphism
$$F':\PreDb(\Ar)\otimes\PreDb(\Ar)\to\PreDb(\bimcat)\fs$$

The codomain $\PreDb(\Ar)\otimes\PreDb(\Ar)$ is not pretriangulated (for example, as in the abelian case, it doesn't contain the cone of $11\to22$) so we take the pretriangulated closure
$$\PreTr(F'):\PreTr(\PreDb(\Ar)\otimes\PreDb(\Ar))\to\PreTr(\PreDb(\bimcat))\fs$$
We also have the obvious map
$$\PreDb(\Ar)\otimes\PreDb(\Ar)\to\PreDb(\Sq)\fs$$
Putting these together we get a commutative diagram
$$\xymatrix{
\PreDb(\Ar)\otimes\PreDb(\Ar)\ar[r]^(0.42){F_1\otimes F_2}\ar[rd]^{F'}\ar[d] &\PreDb(\bimcat)\otimes\PreDb(\bimcat)\ar[d]^{-\otimes_A-}\\
\PreDb(\Sq) & \PreDb(\bimcat)
}$$
of DG-categories which, on applying $\PreTr(-)$ everywhere, gives a commutative diagram
$$\xymatrix{
\PreTr(\PreDb(\Ar)\otimes\PreDb(\Ar)\ar[r]^{}\ar[rd]^{\PreTr(F')}\ar[d]) &\PreTr(\PreDb(\bimcat)\otimes\PreDb(\bimcat)\ar[d]^{})\\
\PreTr(\PreDb(\Sq)) & \PreTr(\PreDb(\bimcat))
}$$
of pretriangulated categories.  On applying $\Ho$ to this diagram, the leftmost vertical map becomes an equivalence by Lemma \ref{tenspre}, so we can define
$$F:\Db(\Sq))\to \Db(E\da E))$$
as the unique map making the diagram
$$\xymatrix{
\Tr(\PreDb(\Ar)\otimes\PreDb(\Ar))\ar[r]^{}\ar[rd]^{\Tr(F')}\ar[d]^\sim &\Tr(\PreDb(E\da E)\otimes\PreDb(E\da E)\ar[d]^{})\\
\Db(\Sq)\ar@{-->}[r]^F & \Db(E\da E))
}$$
commute.  Then by lifting objects and maps of $\Db(\Ar)\otimes\Db(\Ar)$ to $\PreDb(\Ar)\otimes\PreDb(\Ar)$ it is clear that $F$ has the desired properties.
\end{pf}
\end{prop}

\subsection{Products of triangles}
The following statement is modelled on May's axiom (TC3), ``The Braid Axiom for Products of Triangles'' \cite{m}, but we have only recorded the properties that we will use later.  The proof, based on \cite{kn}, makes it clear how to conclude other properties that hold in $\Db(\Sq)$ also hold in $\Db(\bimcat)$.
We point out that our notation is different to May's.
\begin{cor}\label{mayholds}
Given two distinguished triangles
$$X_1\arr{\alpha_X}X_2\arr{\beta_X}X_3\rsar{\gamma_X}$$
$$Y_1\arr{\alpha_Y}Y_2\arr{\beta_Y}Y_3\rsar{\gamma_Y}$$
in $\Db(A\da A)$ we have an object $\kappa\in\Db(A\da A)$ and three more distinguished triangles
$$X_2Y_1\arr{}\kappa\arr{}X_1Y_3\rsar{\alpha_X\gamma_Y}$$
$$\kappa\arr{}X_2Y_2\arr{\beta_X\beta_Y}X_3Y_3\rsar{}$$
$$X_1Y_2\arr{}\kappa\arr{}X_3Y_1\rsar{\gamma_X\alpha_Y}$$
such that the compositions
$$X_2Y_1\to\kappa\to X_2Y_2$$
$$X_1Y_2\to\kappa\to X_2Y_2$$
are $X_2\alpha_Y$ and $\alpha_XY_2$, respectively.

\begin{pf}
Using Lemma \ref{constr}, let $F_X, F_Y:\PreTr(\Ar)\to\PreTr(E\da E)$ send $\alpha:1\to 2$ to $\alpha_X:X_1\to X_2$ and $\alpha_Y:Y_1\to Y_2$ respectively, and then use Proposition \ref{thetrimap} to constuct $F:\Db(\Sq)\to\Db(E\da E)$ which acts as $\Ho(F_1)\otimes\Ho(F_2)$ on the appropriate subcategory.

Applying $\Ho(F_X)$ to the distinguished triangle
$$1\arr{\alpha} 2\arr{\beta} 3\rsar{\gamma} $$
gives a distinguished triangle
$$X_1\arr{\alpha_X}X_2\arr{F_X(\beta)} X'_3\rsar{F_X(\gamma)}$$
which, on comparison with the distinguished triangle
$$X_1\arr{\alpha_X}X_2\arr{\beta_x} X_3\rsar{\gamma_X}$$
tells us that $X'_3:=F_X(3)\cong X_3$.  Similarly, $Y'_3:=F_Y(3)\cong Y_3$.  We conclude that $F(ij)\cong X_iY_j$ for all $i,j\in\{1,2,3\}$.

The distinguished triangles
$$21\arr{} K\arr{} 13\rsar{}$$
$$12\arr{} K\arr{} 31\rsar{}$$
are sent to distinguished triangles
\begin{equation}\label{21k13}
 X_2Y_1\arr{}\kappa\arr{}X_1Y'_3\rsar{}
\end{equation}
\begin{equation}\label{12k31}
 X_1Y_2\arr{}\kappa\arr{}X'_3Y_1\rsar{}
\end{equation}
for some $\kappa\in\Db(A\da A)$.
Then the commutativity of the diagrams
$$
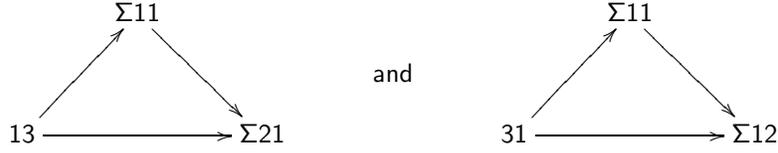
\xymatrix @R=10pt {
 &\shift 11\ar[rdd] &&&& \shift11\ar[rdd] &\\
&&&\text{ and }&&&\\
13\ar[ruu]\ar[rr] &&\shift 21 && 31\ar[ruu]\ar[rr] && \shift12\\
}$$
in $\Db(\Sq)$ implies that their images under $F$ also commute.  From here it is easy to check that the last maps in the distinguished triangles
$$X_2Y_1\arr{}\kappa\arr{}X_1Y_3\rsar{}$$
$$X_1Y_2\arr{}\kappa\arr{}X_3Y_1\rsar{}$$
obtained from distinguished triangles (\ref{21k13}) and (\ref{12k31})
are $\alpha_X\gamma_Y$ and $\gamma_X\alpha_Y$, respectively.

We also have the distinguished triangle
$$K\arr{}22\arr{}33\rsar{}$$
in $\Db(\Sq)$, which is sent by $F$ to a distinguished triangle
$$\kappa\arr{}X_2Y_2\arr{}X'_3Y'_3\rsar{}$$
in $\Db(A\da A)$.  
As the diagrams
$$
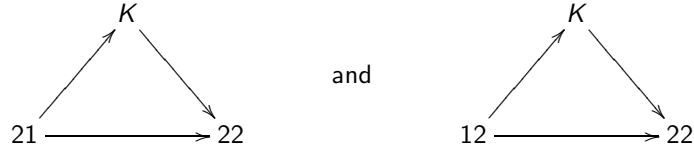
\xymatrix @R=10pt {
 &K\ar[rdd] &&&& K\ar[rdd] &\\
&&&\text{ and }&&&\\
21\ar[ruu]\ar[rr] &&22 && 12\ar[ruu]\ar[rr] &&22\\
}$$
commute we have the final statement.  It only remains to show that the composition $$X_2Y_2\to X'_3Y'_3\arr{\sim} X_3Y_3$$ is $\beta_X\beta_Y$, but this is easy to check given the corresponding statement in $\Db(\Sq)$.
\end{pf}
\end{cor}

Our derived bimodule category also has another nice property:
\begin{cor}\label{mv}
In the above situation, the Mayer-Vietoris triangle
$$X_1Y_1\to X_1Y_2\oplus X_2Y_1\arr{}\kappa\rsar{}$$
is distinguished.
\end{cor}
For a more detailed analysis of such triangles, the interested reader should consult \cite{ive}.

\section{A lifting theorem}\label{sec-lift}
The aim of this section is to show that we can check that relations between periodic twists hold in a certain endomorphism algebra and conclude that they hold more generally.  We will apply this result to braid relations in the following section.

\subsection{Periodic twists}
First we revise the construction of periodic twists from \cite{gra1}.  The construction given here is the same but we will apply it in a more general situation: this greater generality will be important in our proofs.

Recall that if we have an $A\da A$-bimodule $M$ and an algebra automorphism $\tau$ of $A$, then the twisted module $M_\tau$ is defined as the module with the same underlying vector space and left action as $M$, but with right action given by $m\cdot a:=m\tau(a)$ for $m\in M$ and $a\in A$.  Similarly, we can define the module ${}_\tau M$ which is twisted on the left instead of on the right.

Let $P$ be a projective $A$-module and suppose that $E=\End_A(P)^\op$, so $P$ is an $A\da E$-bimodule.  Recall that $E$ is (twisted) periodic if there is a bounded complex $Y$ of projective $E\da E$-bimodules concentrated in degrees $0$ to $n-1$, an algebra automorphism $\tau$ of $E$, and a short exact sequence
$$0\to E_\tau[n-1]\into Y\onto E\to 0$$
of chain complexes of $E\da E$-bimodules.  We say that $E$ has period $n$ and that $Y$ is a truncated resolution of $E$.

Now suppose that we have a short exact sequence
$$0\to F[-1]\into Y\onto E\to0$$
in $\Ch(E\da E)$, where $Y\in\perf(E\da E)$ and we impose no restriction on $F$. 
Denote the map $Y\onto E$ by $f$.  Then by applying the functor $P\otimes_E-\otimes_EP^\vee$ we obtain a map in $\Ch(A\da A)$, and we can consider the composition
$$P\otimes_EY\otimes_EP^\vee\arr{P\otimes f\otimes P^\vee} P\otimes_EE\otimes_EP^\vee \arr{\sim}P\otimes_EP^\vee\arr{\ev}A$$
where $P^\vee$ is the $E\da A$-bimodule $\Hom_A(P,A)$ and the last map is an evaluation map.  We denote this composition by $g:P\otimes_EY\otimes_EP^\vee\to A$.  Then the cone $X$ of $g$ is a bounded complex of $A\da A$-bimodules which are projective both on the left and on the right, and so it induces an endofunctor of the derived category of $A$.

\begin{defn}
Given a short exact sequence
$$0\to F[-1]\into Y\stackrel{f}{\onto} E\to0$$
we have a functor
$$\Psi_{P,f}:=X\otimes_A-:\Db(A)\to\Db(A)\fs$$
\end{defn}

The main result of \cite{gra1} states that if $E$ is periodic and $Y$ is a truncated resolution, so $F$ is a shifted twisted copy of $E$, then $\Psi_{P,f}$ is an autoequivalence.  We call such an autoequivalence a \emph{periodic twist} and write it as $\Psi_{P,Y}$, or just $\Psi_{P}$ when our truncated resolution is minimal.  For an arbitrary $F$, $\Psi_{P,f}$ will not be an autoequivalence, but the construction still gives us a well-defined endofunctor.

\subsection{Periodic twists acting on different derived categories}
Let $P=P_1\oplus\ldots\oplus P_\ell$ be a projective $A$-module which is basic, i.e., it has no two isomorphic nonzero direct summands.  Associate the idempotent $e$ to $P$ and $e_i$ to $P_i$.  Let $E_i=\End_A(P_i)^\op\cong e_iAe_i$ and $E=\End_A(P)^\op\cong eAe$.

Let $Q_i=\Hom_A(P,P_i)\cong eAe_i\cong Ee_i$ for $1\leq i\leq\ell$.  Then $Q_i$ is an $E\da E_i$-bimodule and is projective as a left $E$-module.  In fact, up to isomorphism, all indecomposable projective $E$-modules are obtained in this way.

We collect some basic but important properties of our different projective modules.
\begin{lem}\label{piqi}
 $\End_E(Q_i)\cong\End_A(P_i)$ and $P\otimes_EQ_i\cong P_i$. 

\begin{pf}
We see that
$\End_E(Q_i)\cong e_iEe_i\cong e_i(eAe)e_i=e_iAe_i\cong E_i$
and so $\End_E(Q_i)\cong\End_A(P_i)$.  For the second statement, note that 
$Ae\otimes_E Ee_i\cong Ae_i$.
\end{pf}
\end{lem}

Suppose that for each $1\leq i\leq \ell$ we are given a short exact sequence of chain complexes
$$0\to F_i[-1]\into Y_i\stackrel{f_i}{\onto} E_i\to0$$
where $Y_i\in\perf(E_i\da E_i)$.  We therefore have a distinguished triangle
$$Y_i\arr{f_i}E_i\arr{}F_i\rsa$$
in $\Db(E_i\da E_i)$, where we have used $f_i$ to denote both a map in $\Ch(E_i\da E_i)$ and its image in $\Db(E_i\da E_i)$.

As both the $A$-module $P_i$ and the $E$-module $Q_i$ have endomorphism algebra $E_i$, for each $1\leq i\leq\ell$ we obtain two periodic twists
$$\Psi_{P_i}=X_i\otimes_A-:\Db(A)\arr{}\Db(A)$$
and
$$\Psi_{Q_i}=W_i\otimes_E-:\Db(E)\arr{}\Db(E)$$
where $X_i$ and $W_i$ are defined as the cones of
$$P_i\otimes_{E_i} Y_i\otimes_{E_i} P_i^\vee\arr{g_i} A$$
and
$$Q_i\otimes_{E_i} Y_i\otimes_{E_i} Q_i^\vee\arr{g_i'} E$$
respectively.  The following Lemma, whose proof follows from Lemma \ref{piqi} and the definition of $g_i$ and $g_i'$, will be useful later:
\begin{lem}\label{giprime}
The diagram
$$\xymatrix @C=40pt{
P_i\otimes_{E_i} Y_i\otimes_{E_i} P_i^\vee\ar[r]^(0.7){g_i}\ar[d]^{\sim} & A\\
P\otimes_EQ_i\otimes_{E_i} Y_i\otimes_{E_i} Q_i^\vee\otimes_EP^\vee\ar[r]^(0.7){P\otimes g_i'\otimes P^\vee} & P\otimes_EP^\vee\ar[u]_{\ev_P}
}$$
in $\Ch(A\da A)$ commutes.
\end{lem}

The complexes $X_i$ and $W_i$ fit into distinguished triangles
$$P_i\otimes_{E_i} Y_i\otimes_{E_i} P_i^\vee\arr{g_i} A\arr{h_i} X_i\rsar{i_i}$$
in $\Db(A\da A)$ and
$$Q_i\otimes_{E_i} Y_i\otimes_{E_i} Q_i^\vee\arr{g_i'} E\arr{h_i'} W_i\rsar{i_i'}$$
in $\Db(E\da E)$, with maps denoted as labelled.  We have used the symbol $i$ to mean two different things here: in one context it denotes a map and in another it denotes an integer.  But there should be no confusion as the map $i$ will always have a subscript and will never be used as a subscript.

Sometimes, to save space, we will write
$$V_i=P_i\otimes_{E_i} Y_i\otimes_{E_i}P_i^\vee\in\Db(A\da A)$$
and
$$Z_i=Q_i\otimes_{E_i} Y_i\otimes_{E_i}Q_i^\vee\in\Db(E\da E)$$
so that our distinguished triangles look like
$$V_i\arr{g_i} A\arr{h_i} X_i\rsar{i_i}$$
in $\Db(A\da A)$ and
$$Z_i\arr{g_i'} E\arr{h_i'} W_i\rsar{i_i'}$$
in $\Db(E\da E)$.

\subsection{The lifting theorem}

We want to show that, loosely, periodic twists that decompose ``downstairs'' (i.e., on $\Db(E)$) decompose in the same way ``upstairs'' (i.e., on $\Db(A)$).  We carry over the notation from the previous subsection.  

\begin{lem}\label{yij-tri}
There exists an object $Y_{1,2}$ and 
three triangles
$$Y_{1,2}\arr{f_{1,2}}E\arr{h'_1 h'_2} W_1W_2\rsar{}$$
$$Z_1\arr{p'_1}Y_{1,2}\arr{q'_1}W_1Z_2\rsar{i_1'g_2'}$$
$$Z_2\arr{p'_2}Y_{1,2}\arr{q'_2}Z_1 W_2\rsar{g_1'i_2'}$$
in $\Db(E\da E)$ such that $g'_1=f_{1,2}\circ p'_1$ and $g'_2=f_{1,2}\circ p'_2$.  The second and third of these triangles live in $\per(E\da E)$.

\begin{pf}
Apply Corollary \ref{mayholds} to the triangles
$$Z_i\arr{g_i'} E\arr{h_i'} W_i\rsar{i_i'}$$
for $i=1$ and $2$.

We now show the second and third triangles are perfect.  
As $W_1$ and $W_2$ both have the property of being bounded in nonnegative degrees with $E$ being the only module in its underlying chain complex which is not a projective $E\da E$-bimodule, so does their tensor product.  So the cone of $h_1'h_2':E\to W_1W_2$ is perfect and hence so is $Y_{1,2}$.  The other objects in the second and third triangles are clearly perfect.
\end{pf}
\end{lem}

The next proposition should be compared to Proposition 3.3.3 of \cite{gra1}.
\begin{prop}\label{pdnp}
There is an isomorphism of triangles $P\otimes_E\Delta\cong\nabla\otimes_AP$ in $\Db(A\da E)$, where $\Delta$ and $\nabla$ are the triangles
$$Q_i\otimes_{E_i} Y_i\otimes_{E_i} Q_i^\vee\arr{g_i'} E\arr{h_i'} W_i\rsar{i_i'}$$
and
$$P_i\otimes_{E_i} Y_i\otimes_{E_i} P_i^\vee\arr{g_i} A\arr{h_i} X_i\rsar{i_i}$$
in $\Db(E\da E)$ and $\Db(A\da A)$, respectively.

\begin{pf}
We want a commutative square
$$\xymatrix{
P\otimes_EQ_i\otimes_{E_i} Y_i\otimes_{E_i} Q_i^\vee\ar[d]^\sim\ar[r]^(0.67){P\otimes g'_i} &P\otimes_EE\ar[d]^\sim\\
P_i\otimes_{E_i} Y_i\otimes_{E_i} P_i^\vee\otimes_AP\ar[r]^(0.67){g_i\otimes P} &A\otimes_AP
}$$
where the vertical maps are isomorphisms, and then the proposition will follow by the $5$-lemma for triangulated categories.  We define these maps, and show that the square commutes, by considering the following diagram:
$$\xymatrix{
P\otimes_EQ_i\otimes_{E_i} Y_i\otimes_{E_i} Q_i^\vee\ar[d]^\sim\ar[r] & P\otimes_EQ_i\otimes_{E_i}E_i\otimes_{E_i}Q_i^\vee\ar[r]\ar[d]^\sim &P\otimes_EQ_i\otimes_{E_i}Q_i^\vee\ar[r]\ar[d]^\sim &P\otimes_EE\ar[d]^\sim\\
P_i\otimes_{E_i} Y_i\otimes_{E_i} Q_i^\vee\ar[d]^\sim\ar[r] & P_i\otimes_{E_i}E_i\otimes_{E_i}Q_i^\vee\ar[r]\ar[d]^\sim &P_i\otimes_{E_i}Q_i^\vee\ar[d]^\sim &P\ar[d]^\sim\\
P_i\otimes_{E_i} Y_i\otimes_{E_i} P_i^\vee\otimes_AP\ar[r] &P_i\otimes_{E_i}E_i\otimes_{E_i}P_i^\vee \otimes_AP\ar[r] &P_i\otimes_{E_i}P_i^\vee\otimes_A P\ar[r] 
&A\otimes_AP
}$$
The two leftmost squares commute by the naturality of tensoring with an isomorphism and one can easily check in the module category that the remaining squares and the one hexagon commute.
\end{pf}
\end{prop}

\begin{cor}\label{onetritooth}
$P^\vee\otimes_A\nabla\otimes_AP\cong\Delta$ in $\Db(E\da E)$.

\begin{pf}
$P^\vee\otimes_A\nabla\otimes_AP\cong P^\vee\otimes_AP\otimes_E\Delta\cong E\otimes_E\Delta\cong\Delta\fs$
\end{pf}
\end{cor}

\begin{lem}\label{tri-stu}
There are two triangles
$$PY_{1,2}P^\vee\arr{s}A\arr{}X_1X_2\rsar{}$$
and
$$V_1V_2\arr{\left({V_1g_2}\atop{-g_1V_2}\right)} V_1\oplus V_2\arr{(p_1,p_2)} PY_{1,2}P^\vee\rsar{} $$ 
in $\Db(A\da A)$ such that, for $i\in\{1,2\}$, $s\circ p_i=g_i$ and, after making the obvious identifications, $p_i=Pp'_iP^\vee$.

\begin{pf}
We use Corollary \ref{mayholds} again: take the two distinguished triangles
$$V_i\arr{g_i} A\arr{h_i} X_i\rsar{i_i}$$
for $i=1,2$ and we get distinguished triangles
$$\bar{Y}\arr{\bar{s}}AA\arr{h_1 h_2}X_1X_2\rsar{}$$
$$V_1A\arr{\bar{p}_1}\bar{Y}\arr{\bar{q}_1}X_1V_2\rsar{i_1g_2}$$
$$AV_2\arr{\bar{p}_2}\bar{Y}\arr{\bar{q}_2}V_1X_2\rsar{g_1i_2}$$
where we have labelled the map $\bar{Y}\to AA$ as $\bar{s}$.

Consider the distinguished triangle
$$V_1A\arr{\bar{p}_1}\bar{Y}\arr{\bar{q}_1}X_1V_2\rsar{}$$
from above.
The triangle
$$Z_1\arr{p'_1}Y_{1,2}\arr{q'_1}W_1Z_2\rsar{r'_1}$$
from Lemma \ref{yij-tri} lives in $\per(E\da E)$ so we can apply $P\otimes_E-\otimes_EP^\vee$ to it (without deriving any functors)  to obtain a new distinguished triangle in $\Db(A\da A)$.  Then we want to use isomorphisms as follows
$$\xymatrix @C=60pt {
V_1A\ar[r]^{\bar{p}_1}\ar[d]^\sim &\bar{Y}\ar[r]^{\bar{q}_1}\ar@{-->}[d] &X_1V_2\ar[r]^{i_1g_2}\ar[d]^\sim &V_1A[1]\ar[d]^\sim \\
PZ_1P^\vee\ar[r]^{Pp'_1P^\vee} &P_1Y_{1,2}P_1^\vee\ar[r]^{Pq'_1P^\vee} &PW_1Z_2P^\vee\ar[r]^{Pi'_1g'_2P^\vee} &PZ_1P^\vee[1]
}$$
in order to show $\bar{Y}\cong PY_{1,2}P^\vee$, 
This square
lives in the full subcategory of $\Db(A\da A)$ generated by $P\otimes_kP^\vee$, which we will temporarily denote $\Db(P\da P^\vee)$.  The functor
$$P^\vee\otimes_A-\otimes_AP:\Db(A\da A)\to \Db(E\da E)$$
restricts to an equivalence on $\Db(P\da P^\vee)$.  By Corollary \ref{onetritooth} we have that $P^\vee\otimes_Ag_i\otimes_AP=g'_i$ and $P^\vee\otimes_Ai_i\otimes_AP=i'_i$, so the square is sent to a commutative square after applying $P^\vee\otimes_A-\otimes_AP$, hence it must have been commutative to start with.

Now we have a map $\bar{s}:\bar{Y}\to A\otimes_AA$ and an isomorphism $\bar{Y}\cong PY_{1,2}P^\vee$, so we use $s$ to label the composition
$$s:PY_{1,2}P^\vee\to \bar{Y}\to A\otimes_AA\arr{\sim}A\fs$$
Then we define the distinguished triangle
$$PY_{1,2}P^\vee\arr{s}A\arr{}X_1X_2\rsar{}$$
by the commutative diagram
$$\xymatrix @C=60pt {
PY_{1,2}P^\vee\ar[r]^{s}\ar[d]^\sim &A\ar[r]^{}\ar[d]^\sim &X_1X_2\ar@{~>}[r]^{}\ar@{=}[d] & \\
\bar{Y}\ar[r]^{\bar{s}} &AA\ar[r]^{} &X_1X_2\ar@{~>}[r]^{} &
}$$

We also get the Mayer-Vietoris triangle
$$V_1V_2\arr{\left({V_1g_2}\atop{-g_1V_2}\right)} V_1\oplus V_2\arr{(p_1,p_2)} PY_{1,2}P^\vee\rsar{} $$ 
from Corollary \ref{mv}, where $p_i$ is the composition $$V_i\arr{\sim}V_iA\arr{\bar{p_i}}\bar{Y}\arr{\sim}PY_{1,2}P^\vee$$ and similarly for $q_i$.

One can check that the equations $s\circ p_i=g_i$ follow from $\bar{s}\circ\bar{p_1}=g_1\otimes_AA$ and $\bar{s}\circ\bar{p_2}=A\otimes_Ag_2$, and it is clear from our earlier isomorphism of triangles that $p_i=Pp_i'P^\vee$, i.e., the diagram
$$\xymatrix{
V_i \ar[r]^{p_i}\ar[d]^{\sim} & PY_{1,2}P^\vee \\
V_iA\ar[r]^{\bar{p_i}}\ar[d]^{\sim} & {} \bar{Y}\ar[u]_{\sim}\ar[d]^{\sim} \\
PZ_1P^\vee\ar[r]^{Pp'_1P^\vee} &P_1Y_{1,2}P_1^\vee
}$$
commutes.
\end{pf}
\end{lem}

\begin{prop}
$\Psi_{P_1,f_1}\circ\Psi_{P_2,f_2} =\Psi_{P,f_{1,2}}\fs$

\begin{pf}
Let $X_{1,2}$ be the cone of the composition $PY_{1,2}P^\vee\arr{Pf_{1,2}P^\vee}PP^\vee\arr{\ev}A$.  Then we want to show that $X_{1,2}\cong X_1X_2$.
As usual, we prove this by constructing a commutative diagram
$$\xymatrix{
PY_{1,2}P^\vee\ar[rr]^{s}\ar@{=}[d] & & A\ar@{=}[d]\\
PY_{1,2}P^\vee\ar[rd]^{Pf_{1,2}P^\vee} & & A\\
 & PP^\vee\ar[ur]^{\ev_P} &
}$$
and appealing to the triangulated $5$-lemma, so we need to show that the map
$$v=s-Pf_{1,2}P^\vee\circ\ev_P:PY_{1,2}P^\vee\to A$$
is zero.  We argue more indirectly this time.

By the previous lemma, we have 
a distinguished triangle
$$\xymatrix @R=6pt {
 & V_1\ar[rd]^{p_1} & & \\
V_1V_2\ar[ru]^{V_1g_2}\ar[rd]_{-g_1V_2} & \oplus & PY_{1,2}P^\vee\ar@{~>}[r] & \\
 & V_2\ar[ru]_{p_2} & & 
}$$
in $\Db(A\da A)$. 
Suppose we can show that $(p_1,p_2)\circ v=0$.  Then the commutative diagram
$$\xymatrix @R=6pt {
 V_1\ar[rd]^{p_1}\ar@/^1.4pc/[dddddd] |!{[dd];[rr]}\hole & & & \\
 \oplus & PY_{1,2}P^\vee\ar[r]\ar[ddddd]^{v} &V_1V_2[1]\ar@{~>}[r]\ar@{-->}[ddddd] & \\
 V_2\ar[ru]_{p_2}\ar[dddd] & & & \\
 & & & \\
 & & & \\
 & & & \\
0\ar[r] & A\ar@{=}[r] & A\ar@{~>}[r] &
}$$
and the completion axiom for triangulated categories shows that $v$ must factor through $V_1V_2[1]$.  But $A$ is concentrated in degree $0$, and $V_1V_2[1]$ is concentrated in strictly negative degrees, so $v$ must be zero.

It remains to show that $(p_1,p_2)\circ v=0$, or, equivalently, that
for $i=1,2$, $$p_i\circ s=p_i\circ Pf_{1,2}P^\vee\circ \ev_P\fs$$
By the previous lemma, the left hand side is $g_i$ and the right hand side is $P(p_i'\circ f_{1,2})P^\vee\circ \ev_P$, which by Lemma \ref{yij-tri} is equal to $Pg_i'P^\vee\circ \ev_P$.  But this is equal to $g_i$, as Lemma \ref{giprime} tells us that the corresponding statement is true in the chain complex category.
\end{pf}
\end{prop}

If $\ii=(i_1,\ldots,i_r)\in\{1,\ldots,\ell\}^r$ then we will write $\Psi_\ii=
\Psi_{P{i_1}}\Psi_{P{i_2}}\ldots \Psi_{P{i_r}}:\Db(A)\arr{}\Db(A)$ and similarly for $\Psi_\ii'
:\Db(E)\arr{}\Db(E)$.  Recall that $P=P_1\oplus\cdots\oplus P_\ell$ and let $Q=Q_1\oplus\cdots\oplus Q_\ell$.

\begin{thm}[Lifting theorem]\label{lifting}
Suppose we have a short exact sequence
$$0\to F[-1]\into Y\stackrel{f}{\onto} E\to0$$
in $\Ch(E\da E)$, with $Y\in\per(E\da E)$, and we have a natural isomorphism
$\Psi'_\ii\cong F\dert_E-$
of functors.  Then
$$\Psi_\ii\cong \Psi_{P,f}\fs$$
In particular,
\begin{enumerate}
\item if $\Psi_\ii'\cong\Psi_\jj'$ for $\jj=(i_1,\ldots,i_r)\in\{1,\ldots,\ell\}^r$ then  $\Psi_\ii\cong\Psi_\jj$, and
 \item if $\Psi_\ii'\cong E_\sigma[n]\otimes_E-$, so $E$ is twisted-periodic, then $\Psi_\ii$ is isomorphic to the periodic twist $\Psi_P$.
\end{enumerate}

\begin{pf}
This follows from the previous proposition by induction.
\end{pf}
\end{thm}

\section{Braid relations and longest elements of symmetric groups}\label{sec-long}
Suppse we have a braid group acting by spherical twists on the derived category of a symmetric algebra.  Then using the lifting theorem developed in the previous section, we will show that lifts of longest elements from symmetric groups to braid groups act in the way suggested by the example at the end of \cite{gra1}.

\subsection{Brauer tree algebras of lines without multiplicity}\label{gamman}

We define a collection of algebras $\Gamma_n$, $n\geq1$, as path algebras of quivers with relations.  Let $\Gamma_1=k[x]/\gen{x^2}$ and let $\Gamma_2=kQ_2/I_2$, where $Q_2$ is the quiver
$$\xymatrix@=10pt{
  Q_2&=&1\ar@/^/[rr]^{\alpha} & & 2\ar@/^/[ll]^\beta
}$$
and $I_2$ is the ideal generated by $\alpha\beta\alpha$ and $\beta\alpha\beta$.  For $n\geq3$, let $\Gamma_n=kQ_n/I_n$ where $Q_n$ is the quiver
$$\xymatrix@=10pt{
  Q_n&=&1\ar@/^/[rr]^{\alpha_1} & & 2\ar@/^/[ll]^{\beta_2}\ar@/^/[rr]^{\alpha_2} & & {}\cdots\ar@/^/[ll]^{\beta_3}\ar@/^/[rr]^{\alpha_{n-1}} & & n\ar@/^/[ll]^{\beta_n}
}$$
and $I_n$ is the ideal generated by $\alpha_{i-1}\alpha_{i}$, $\beta_{i+1}\beta_i$, and $\alpha_i\beta_{i+1}-\beta_i\alpha_{i-1}$ for $2\leq i\leq n-1$.  
For ease of notation, let $\alpha=\sum\alpha_i$ and $\beta=\sum\beta_i$. 
Then we can write $\alpha_i$ and $\beta_j$ as $e_i\alpha$ or $\alpha e_{i+1}$ and $e_j\beta$ or $\beta e_{j-1}$.

The algebras $\Gamma_n$ have appeared in many contexts.  Some examples are:
\begin{itemize}
\item $\Gamma_n$ is the trivial extension algebra of the path algebra of a Dynkin quiver of type $A_n$ with bipartite orientation \cite{bbk,hk};
\item $\Gamma_n$ is quadratic dual to the preprojective algebra of type $A_n$ with the path-length grading, for $n>2$ \cite{bbk,hk};
\item $\Gamma_n$ is a Brauer tree algebra of a line
$$\xymatrix{\circ\ar@{-}[r] &\circ\ar@{-}[r] &\cdots\ar@{-}[r] &\circ\ar@{-}[r] &\circ}$$
 with $n$ edges and no exceptional vertex; up to derived equivalence, these are all Brauer tree algebras with multiplicity function the constant $1$ \cite{rz};
\item $\Gamma_n$ is the zig-zag algebra of type $A_n$, which was used by Huerfano and Khovanov to categorify the adjoint representation of a type $A$ quantum group \cite{hk};
\item $\Gamma_n$ is isomorphic to the underlying ungraded algebra of the formal differential graded ext-algebras of an $A_n$-configurations of spherical objects \cite{st}.
\end{itemize}

The relations for $\Gamma_n$ are homogeneous, and there are various possible ways to put a grading on these algebras (see, for example, \cite[Section 4]{st}).  When we want to consider them as graded algebras, we will give $x$ in $\Gamma_1$ degree $2$, and for $n\geq2$, all $\alpha_i$ and $\beta_j$ will have degree $1$, as in \cite{hk}.

Let $e_i$ denote the primitive idempotent corresponding to the vertex $i$ of $Q_n$.  We have an algebra automorphism $\tau_n\in\Aut(\Gamma_n)$ of order $2$: $\tau_1$ sends $x$ to $-x$, $\tau_2$ swaps $e_i$ and $e_{n+1-i}$ and $\alpha$ and $\beta$,
 and for $n>2$, $\tau_n$ sends the idempotent $e_i$ to $e_{n+1-i}$ and swaps $\alpha_i$ and $\beta_{n+1-i}$. 
 Note that these automorphisms respect our grading.

Let $P_i$ denote the projective $\Gamma_n$-module $\Gamma_ne_i$.  With our grading conventions the following result, which is easy to prove, is true in both the graded and ungraded setting:
\begin{lem}\label{gammaend}
Let $1\leq i\leq j\leq n$.  Then
$$\End_{\Gamma_n}(P_i\oplus P_{i+1}\oplus\ldots\oplus P_j)^\op\cong\Gamma_{j-i+1}\fs$$
In particular, for $1\leq i\leq n$, $$\End_{\Gamma_n}(P_i)\cong k[x]/\gen{x^2}\fs$$
\end{lem}

\subsection{Spherical twists for symmetric algebras}
Let $A$ be a symmetric algebra.
\begin{defn}
We say that a projective $A$-module is \emph{spherical} if $\End_A(P)\cong k[x]/\gen{x^2}$.
\end{defn}
By Lemma \ref{gammaend}, for $n\geq1$, each indecomposable projective $\Gamma_n$-module is spherical.

\begin{thm}[\cite{rz} for the algebras $\Gamma_n$; \cite{st} in general]
If $P$ is a projecive $A$-module 
which is spherical
then the functor $F_P$ given by tensoring with the complex
$$P\otimes_k P^\vee\arr{\ev} A$$
of $A\da A$-bimodules concentrated in degrees $1$ and $0$ is a derived autoequivalence.
\end{thm}
These equivalences are called spherical twists.  Note that they are special cases of periodic twists.

Recall that the braid group $B_{n+1}$ on $n+1$ letters has the presentation
$$B_{n+1}=\left<\:\sigma_1,\ldots,\sigma_n\:|\:\sigma_i\sigma_j=\sigma_j\sigma_i\text{ for }\abs{i-j}>1;\:\sigma_i\sigma_{i+1}\sigma_i=\sigma_{i+1}\sigma_i\sigma_{i+1}\text{ for }1\leq i<n\; \right>\fs$$

\begin{defn}[\cite{st}]
We say that a collection $\{P_1,\ldots,P_n\}$ of projective $A$-modules is an \emph{$A_n$-configuration} if each $P_i$ is spherical and, for all $1\leq i,j\leq n$,
$$
\dim_k\Hom_A(P_i,P_j) =
\begin{cases}
1 & \text{if } \abs{i-j}=1\semic\\
0 & \text{if } \abs{i-j}>1\fs
\end{cases}$$
\end{defn}
We have the following observation, which is straightforward in our setting, and is a special case of the more general statement \cite[Lemma 4.10]{st}:
\begin{lem}\label{anconfiggamma}
Let $A$ be a symmetric algebra.  A collection $\{P_1,\ldots,P_n\}$ of projective $A$-modules is an $A_n$-configuration if and only if
\[\End_A(\bigoplus_{i=1}^nP_i)^\op\cong\Gamma_n\fs\]
\end{lem}

\begin{thm}[\cite{rz} for the algebras $\Gamma_n$; \cite{st} in general]\label{bga}
If the collection $\{P_1,\ldots,P_n\}$ of projective $A$-modules is an $A_n$-configuration then we have an action of the braid group on the derived category of $A$,
$$B_{n+1}\to\Aut(\Db(A))\cma$$
which sends the braid group generator $\sigma_i$ to the spherical twist $F_i$ associated to the projective $A$-module $P_i$.
\end{thm}

As a corollary of the lifting theorem \ref{lifting}, we obtain a new proof that the spherical twists associated to an $A_n$-configuration satisfy the braid relations.

\emph{Proof of Theorem \ref{bga}:}
It is easy to see from the definitions that if $\abs{i-j}>1$ then $F_iF_j\cong F_jF_i$; the hard part is to show that $F_iF_{i+1}F_i\cong F_{i+1}F_iF_{i+1}$ for $1\leq i\leq n-1$.  But can check directly that $F_1F_2F_1\cong F_2F_1F_2$ on the derived category of $\Gamma_2\cong \End_A(P_i\oplus P_{i+1})$ and, by the lifting theorem, this is enough. \endpf

The algebras $\Gamma_n$ are of finite representation type, and hence are twisted periodic, but in fact we can say more.
\begin{thm}[\cite{bbk}]\label{bbkres}
The algebra $\Gamma_n$ is twisted periodic with period $n$ and automorphism $\tau_n$.
\end{thm}
A natural question is: what do the associated periodic twists look like?  It was noted in \cite{gra1} that periodic twists on $\Gamma_3$ associated to the direct sum of the first two projectives, which have endomorphism algebra $\Gamma_2$, are isomorphic to the composition $F_1F_2F_1$ of spherical twists.  We will show that this pattern continues.

\subsection{Longest elements}
The symmetric group $S_{n+1}$ is the group of automorphisms of the set $\{1,2,\ldots,n+1\}$.  It has a standard generating set consisting of the tranpositions $s_i$ which interchange the numbers $i$ and $i+1$.  Note that these generators are involutions.  There is a minimal number of letters from the alphabet $\{s_i\}_{i=1}^n$ needed to write a given element $w$ of $S_{n+1}$: this number is called the length of $w$.  There is a unique element of $S_{n+1}$ of longest length, called the longest element.  It sends $i\in\{1,2,\ldots,n+1\}$ to $n+2-i$.  We write it as $w_0$ or, when we need to make explicit the dependence on $n$, as $w_0^{(n+1)}$.

A reduced expression for $w\in S_{n+1}$ is a way to write $w$ using the minimal number of standard generators.  There are different reduced expressions for the longest element.  We give one inductive way to write the reduced expression: let $w_0^{(1)}$ be the identity element of the identity group $S_1$, and then
$$w_0^{(n+1)}=w_0^{(n)}s_n\ldots s_2s_1$$
is the longest element of $S_{n+1}$.  Here, we have used the embedding $S_n\into S_{n+1}$ which sends $s_i\in S_n$ to $s_i\in S_{i+1}$ to consider $w_0^{(n)}$ as an element of $S_{n+1}$, and we are employing a slight abuse of notation in using $w_0^{(n)}$ to represent both an element of the group $S_n$ and a word in the alphabet $\{s_i\}_{i=1}^{n-1}$.

The symmetric group has a presentation in terms of the generators $\{s_i\}_{i=1}^n$ as follows:
$$S_{n+1}=\left<\:s_1,\ldots,s_n\:|\:s_is_j=s_js_i\text{ for }\abs{i-j}>1;\:s_is_{i+1}s_i=s_{i+1}s_is_{i+1}\text{ for }1\leq i<n;\:s_i^2=1\; \right>\fs$$
There is an obvious surjective group homomorphism $B_{n+1}\onto S_{n+1}$ from the braid group to the symmetric group which sends $\sigma_i$ to $s_i$.  By sending $s_i$ to $\sigma_i$ we obtain a positive lift of each word in the alphabet $\{s_i\}_{i=1}^n$.  By Matsumoto's Theorem \cite{matsumoto}, we can lift each element of $S_{n+1}$ to $B_{n+1}$ using reduced expressions, and this is independent of our particular choice of reduced expression. 
 We denote the positive lift of $w_0^{(n+1)}\in S_{n+1}$ by $t_{n+1}$.

Consider the algebras $\Gamma_n$ described in Subsection \ref{gamman}.  By Lemma \ref{gammaend}, each projective $\Gamma_n$-module $P_i$ has endomorphism algebra $k\gen{x}/\gen{x^2}$, and so we have a spherical twist given by the period $1$ twisted resolution of $E_i=\End_{\Gamma_n}(P_i)^\op$, which we label $F_i:\Db(\Gamma_n)\arr{\sim}\Db(\Gamma_n)$.  By Lemma \ref{anconfiggamma} the indecomposable projective $\Gamma_n$-modules form an $A_n$-configuration and so we have a braid group action
$$\varphi_n:B_{n+1}\to \Aut(\Db(\Gamma_n))$$
$$\sigma_i\mapsto F_i\fs$$
It is natural to ask what the image of $t_{n+1}$ is under $\varphi_n$.  The answer was given by Rouquier and Zimmermann:
\begin{thm}\cite[Theorem 4.5]{rz}\label{rzthm}
 The longest element $t_{n+1}$ of $B_{n+1}$ acts on $\Db(\Gamma_n)$ as the shift and twist $-_{\tau_n}[n]$.
\end{thm}

Note that a differential graded analogue of this theorem appears as Lemma 3.1 in \cite{sei}.  In the diffential graded setting, many of the technical difficulties in proving this theorem disappear.

We can combine the lifting theorem \ref{lifting} and Theorem \ref{rzthm} to answer the more general question for an arbitrary symmetric algebra $A$: given an $A_n$-configuration, how does the longest element $t_{n+1}$ act on $\Db(A)$?
\begin{cor}
For a symmetric algebra $A$ and an $A_n$-configuration $\{P_1,\ldots,P_n\}$, the longest element $t_{n+1}$ acts as the periodic twist $\Psi_P$, where $P=P_1\oplus\cdots\oplus P_n$.
\end{cor}

\subsection{Quadratic algebras and Koszul complexes}\label{ss-quad}
This section recounts well-knows ideas and constructions from Koszul duality (\cite{pri}, \cite{bgs}, \cite{bg}, \cite{bking}, \ldots), but we need to work with arbitrary quadratic algebras instead of only Koszul algebras (see \cite{mos}).  Our aim is to use this theory, as well as the theory of almost Koszul duality due to Brenner-Butler-King \cite{bbk}, to further study the algebras $\Gamma_n$.

Let $\Lambda=\bigoplus_{i\in\Z}\Lambda_i$ be a $\Z$-graded $k$-algebra.  We will assume that $\Lambda$ is positively graded and generated in degree $1$, i.e., $\Lambda=\bigoplus_{i\geq0}\Lambda_i$ and for $i>0$, $\Lambda_i=\Lambda_1\Lambda_{i-1}=\Lambda_{i-1}\Lambda_1$.  A $\Lambda$-module $M=\bigoplus_{i\in\Z}M_i$ is graded if $\Lambda_iM_j\subseteq M_{i+j}$ for all $i,j$.  Let $\Lambda\grmod$ denote the category of finitely generated graded $\Lambda$-modules, in which all maps $f:M\to N$ are homogeneous of degree $0$, i.e., $f(M_i)\subseteq N_i$.  For each $n\in\Z$ there is an autoequivalence $\grsh n$ of $\Lambda\grmod$ which sends $M$ to the module $M\grsh n$ with $M\grsh{n}_i=M_{n+i}$.  If each $M_i$ is finite dimensional over $k$, then we write $M^*$ for $\bigoplus_{i\in\Z} (M_i)^*$.  Note that $(M^*)_i=(M_{-i})^*$.

\begin{defn}
 We say that $\Lambda$ is \emph{Frobenius of Gorenstein parameter $n$} if we have an isomorphism $\varphi:\Lambda\arr{\sim}{}_{\nu}(\Lambda^*)\grsh{-n}$ in the category $\Lambda\grmodgr\Lambda$ of finitely generated graded $\Lambda\da\Lambda$-bimodules for some graded algebra automorphism $\nu\in\Aut(\Lambda)$, called the \emph{Nakayama automorphism}.  If moreover $\nu$ is the identity, we say that $\Lambda$ is \emph{symmetric of Gorenstein parameter $n$}.
\end{defn}

Let $S$ be a semisimple $k$-algebra $S$ which, for simplicity, we will assume is basic, i.e., it is a product of copies of the field $k$.  Let $V$ be an $S\da S$-bimodule.  Then recall that the tensor algebra $\Tens_S(V)$ is a positively graded algebra
$$\Tens_S(V)=\bigoplus_{i\geq0}V^{\otimes_S i}\cma$$
where $V^{\otimes_S 0}=k$ and $V^{\otimes_S i}=V\otimes_S V\otimes_S \cdots\otimes_S V$ with $i$ factors.  The multiplication in $\Tens_S(V)$ is given by the obvious concatenation.

\begin{defn}\label{defnquad}
We say $\Lambda$ is \emph{quadratic} if there is a semisimple $k$-algebra $S$, an $S\da S$-bimodule $V$, and a subset $R\subset V\otimes_SV$ such that $\Lambda\cong \Tens_S(V)/(R)$.
\end{defn}
Note that as $(R)$ is a homogeneous ideal, quadratic algebras inherit a positive grading from the tensor algebra.

For the rest of this subsection, assume $\Lambda$ is quadratic.  The quadratic dual $\Lambda^!$ is defined as
$$\Lambda^!=\Tens_S(V^*)/(R^\perp)$$
where $R^\perp=\{f\in (V\otimes_S V)^*\:|\:f(R)=0\}$ is the perpendicular space to $R$ and, for $V,W\in S\grmodgr S$, we identify the $S\da S$-bimodules $(V\otimes_S W)^*$ and $W^*\otimes_S V^*$.
$\Lambda^!$ is also a quadratic algebra, and $(\Lambda^!)^!\cong\Lambda$.

If we have a graded algebra automorphism $\tau:\Lambda\arr{\sim}\Lambda$ then restricting to degrees $0$ and $1$ gives an algebra automorphism $\tau_0:S\arr{\sim} S$ and a vector space automorphism $\tau_1:V\arr{\sim} V$. 
We have a $k$-algebra map $\tau^!_0:=\tau_0:S\arr{\sim} S$ and taking the dual of the inverse of $\tau_1$ gives us a vector space map $\tau_1^!:V^*\arr{\sim} V^*$.  These maps generate an algebra automorphism $\tau^!:\Lambda^!\arr{\sim}\Lambda^!$, and we have $(\tau^!)^!=\tau$.

Let $M=\bigoplus M_i\in \Lambda^!\grmodgr \Lambda^!$ be a graded $\Lambda^!\da\Lambda^!$-bimodule, so we have left and right actions $\Lambda^!\otimes_SM\to M$ and $M\otimes_S\Lambda^!\to M$.  Restricting to $(\Lambda^!)_1=V^*$ gives maps $\ell_M:V^*\otimes_S M_{i-1}\to M_{i}$ and $r_M:M_{i-1}\otimes_SV^* \to M_{i}$ which have duals
$\ell_M^*:(M_i)^*\to (M_{i-1})^*\otimes_S V$
and
$r_M^*:(M_i)^*\to V\otimes_S (M_{i-1})^*$, where we have used the natural isomorphism $V^{**}\cong V$ to identify $V$ and its double dual.

We say that $M\in \Lambda\grmod$ is generated in degree $i$ if there exists some subset $L\subset M_i$ such that $M=\Lambda L\Lambda$.
Let $\Lambda\grprojj\Lambda$ denote the additive category of graded projective $\Lambda\da\Lambda$-bimodules  and let $\lin(\Lambda\grprojj\Lambda)$ denote the corresponding category of linear complexes: this is the category of chain complexes
$$\cdots\to X_{1}\to X_0\to X_{-1}\to \cdots$$
where $X_i\in\Lambda\grprojj\Lambda$ is generated in degree $i$ and all maps are homogeneous of degree $0$.

\begin{defn}
There is a contravariant functor
$$Q:\Lambda^!\grmodgr \Lambda^!\to\lin(\Lambda\grprojj \Lambda)$$
which is defined as follows:
$Q$ sends $M=\bigoplus_{i\in\Z}M_i\in\Lambda^!\grmodgr \Lambda^!$ to the complex
$$\cdots \arr{} Q(M)_1\arr{d}Q(M)_0 \arr{d}Q(M)_{-1}\to\cdots $$
with
$ Q(M)_i= \Lambda\otimes_S (M_i)^*\otimes_S\Lambda\grsh{-i}$ (where we consider $M_i$ as an $S\da S$-bimodule concentrated in grade $0$) and
differential $d$ given by the following composition:
$$\xymatrix{
 &\Lambda(M_{i-1})^* V \Lambda\grsh{-i+1}\ar[rd]^{\Lambda(M_{i-1})^*  m} &\\
\Lambda(M_i)^* \Lambda\grsh{-i}\ar[ur]^{\Lambda \ell_M^* \Lambda}\ar[dr]^{(-1)^i\Lambda r_M^* \Lambda} && \Lambda(M_{i-1})^* \Lambda\grsh{-i+1} \\
 &\Lambda V (M_{i-1})^* \Lambda\grsh{-i+1}\ar[ru]^{m(M_{i-1})^* \Lambda} &
}$$
Here, $m:\Lambda\otimes_S\Lambda\to \Lambda$ is the algebra multiplication map and in applying $m$ we have implicitly used the inclusion $V\into \Lambda$.
\end{defn}
This construction does give a chain complex because the differential is defined so that, of the four ways an element can be mapped from $P_{i+1}\grsh{i+1}$ to $P_{i-1}\grsh{i-1}$, two of these are zero by definition of the quadratic dual and the other two cancel out due to the choice of sign.  The contravariance of $Q$ comes from the contravariance of $(-)^*$ and the functoriality comes from the fact that bimodule maps give commuting chain complex maps.

In the case when $Q$ is Koszul, the functor $Q$ can be extended to complexes of $\Lambda^!\da\Lambda^!$-bimodules and we recover the Koszul duality functor of \cite[Theorem 2.12.1]{bgs}.

Note that $\Lambda^!\grmodgr \Lambda^!$ and $\lin(\Lambda\grprojj \Lambda)$ are both abelian categories.
\begin{lem}\label{Qexact}
$Q$ is an exact functor, i.e., a short exact sequence
$$0\to K\into L\onto M\to 0$$
in $\Lambda^!\grmodgr \Lambda^!$ is sent to a short exact sequence
$$0\to Q(M)\into Q(L)\onto Q(K)\to0$$
in 
$\lin(\Lambda\grprojj \Lambda)$.

\begin{pf}
This follows from the exactness of tensoring over the semisimple algebra $S$.
\end{pf}
\end{lem}

$Q$ plays nicely with various operations we can perform on bimodules.
\begin{prop}\label{Qprop}
Let $M\in\Lambda^!\grmodgr \Lambda^!$.  The functor $Q$ has the following properties:
\begin{enumerate}
 \item\label{Qshift} For $i\in \Z$, $Q(M\grsh{i})=Q(M)\grsh{i}[-i]$;
 \item\label{Qtwist} Let $\tau^!\in\Aut(\Lambda^!)$.  Then $Q(M_{\tau^!})\cong{}_{\tau}Q(M)$.
 \item\label{Qdual} If $\Lambda$ is a graded Frobenius algebra of Gorenstein parameter $n$ then there is a natural isomorphism $Q(M^*)\cong {}_\nu (Q(M)^*)_{\nu^{-1}}\grsh{-2n}$.  In particular, if $\Lambda$ is symmetric then $Q(M^*)\cong Q(M)^*\grsh{-2n}$.
\end{enumerate}

\begin{pf}
\begin{itemize}
 \item[\ref{Qshift}] This is easy to check, and was noted in \cite[Theorem 2.12.5(ii)]{bgs}.
 \item[\ref{Qtwist}] For notational simplicity, let $N=M_{\tau^!}$.  Then we have an isomorphism of $\Lambda\da\Lambda$-bimodules
$$\tau\otimes_S\id\otimes_S\id:\Lambda\otimes_S (N_i)^*\otimes_S\Lambda\arr{\sim}{}_\tau \Lambda\otimes_S (M_i)^*\otimes_S\Lambda\fs$$
To see that this commutes with the differential, use the fact that $\tau$ is an algebra homomorphism and note that $r_N=r_M\circ (\id\otimes \tau^!_1)$, so $r_N^*=(\tau_1^{-1}\otimes\id)\circ r_M^*$.
 \item[\ref{Qdual}] First, assume that $\Lambda$ is graded symmetric, i.e., we are given an isomorphism $\varphi:\Lambda\arr{\sim}\Lambda^*\grsh{-n}$ of $\Lambda\da\Lambda$-bimodules. 
We have a isomorphisms
$$\varphi\otimes\id\varphi^*\grsh{-n}:Q(M^*)_i\cong (Q(M)^*)_i\grsh{-2n}$$
for each $i\in\Z$,
so we just need to show that the square
$$\xymatrix{
\Lambda(M_{i})^*\Lambda\ar[r]^d\ar[d]^{\varphi.1.\varphi^*\grsh{-n}}\grsh{i} & \Lambda(M_{i+1})^*\Lambda\grsh{i+1}\ar[d]^{\varphi.1.\varphi^*\grsh{-n}}\\
\Lambda^*(M_{i})^*\Lambda^*\ar[r]^{d^*}\grsh{i-2n} & \Lambda^*(M_{i+1})^*\Lambda^*\grsh{i+1-2n}
}$$
commutes.  All our maps will be homogeneous of degree $0$, so we will drop the gragings from the notation.

Recall that $d=d^\ell + (-1)^id^r:Q(M)_i\to Q(M)_{i-1}$, where $d^\ell=(1.1.m)\circ (1.\ell_M^*.1)$.  We will show that $(\varphi.1.\varphi^*)\circ d^\ell=(d^r)^*\circ (\varphi.1.\varphi^*)$; the corresponding statement with $d^\ell$ and $d^r$ interchanged is similar.

Let's draw the diagram for  $(\varphi.1.\varphi^*)\circ d^\ell=(d^r)^*\circ (\varphi.1.\varphi^*)$:
$$\xymatrix{
\Lambda(M_{-i})\Lambda \ar[d]^{\varphi.1.\varphi^*}\ar[r]^{1.\ell_{M^*}^*.1} & \Lambda(M_{-i+1})V\Lambda\ar[r]^{1.1.m} & \Lambda(M_{-i+1})\Lambda\ar[d]^{\varphi.1.\varphi^*}\\
\Lambda^*(M_{-i})\Lambda^*\ar[r]^{1.1.m^*} & \Lambda^*(M_{-i})V^*\Lambda^*\ar[r]^{1.r_M.1} &\Lambda^*(M_{-i+1})\Lambda^*
}$$
By definition, we have
$$\ell_{M^*}= (\ev_V.1) \circ(1.r_M^*):V^*\otimes_S(M^*)_{i}\to(M^*)_{i+1}$$
and it is clear that the diagram commutes in the first tensor factor, so we need to show that
$$\xymatrix{
(M_{-i})\Lambda \ar[rr]^{1.\ev^*.1}\ar[dd]^{1.\varphi^*} && (M_{-i})V^*V\Lambda\ar[rr]^{r_M.1.1}\ar[rd]^{1.1.m} &&(M_{1-i})V\Lambda\ar[rd]^{1.m} & \\
&&& (M_{-i})V^*\Lambda\ar[rr]^{r_m.1}\ar[ld]^{1.1.\varphi^*} && (M_{1-i})\Lambda\ar[ld]^{1.\varphi^*} \\
(M_{-i})\Lambda^*\ar[rr]^{1.m^*} && (M_{-i})V^*\Lambda^*\ar[rr]^{r_M.1} && (M_{1-i})\Lambda^* &
}$$
commutes.  The two squares on the right obviously commute so it only remains to check the pentagon on the left.  But, after removing the $M_{-i}$ tensor factor on the left, this is exactly dual to the commutative diagram which states that $\varphi$ is a map of right $A$-modules.

Finally, if $\Lambda$ is graded Frobenius but not symmetric, then we have an isomorphism $\varphi:\Lambda\arr{\sim}{}_\nu(\Lambda^*)$ of $\Lambda\da\Lambda$-bimodules, where $\nu\in\Aut(\Lambda)$ is the Nakayama automorphism.  So $\varphi^*:\Lambda_\nu\arr{\sim}\Lambda^*$, and we use the maps $\varphi.1.(\varphi^*)_{\nu^{-1}}$ as above to show $Q(M^*)\cong {}_\nu (Q(M)^*)_{\nu^{-1}}$.
\end{itemize}
\end{pf}
\end{prop}

\subsection{Quotients and idempotents}
For an algebra $\Lambda$ and an idempotent $e\in\Lambda$, we can define two algebras: $\Lambda/\Lambda(1-e)\Lambda$ and $e\Lambda e$. If $\Lambda$ is semisimple then they are isomorphic but in general they are different.
The first is, by definition, a quotient of $\Lambda$ while the second is, in general, not.  We investigate the relation when $\Lambda$ is quadratic.

\begin{lem}
  Let $\Lambda$ be a quadratic algebra and $e\in\Lambda_0\subset\Lambda$ be an idempodent.  Then $\Lambda/\Lambda(1-e)\Lambda$ is quadratic.  Moreover, if the algebra $e\Lambda e$ is generated in degree $1$ and is quadratic then we have an isomorphism of graded algebras
$$e\Lambda e\cong \left(\frac{\Lambda^!}{\Lambda^!(1-e)\Lambda^!}\right)^!\fs$$

\begin{pf}
 Let $\Lambda\cong \Tens_S(V)/(R)$ as in Definition \ref{defnquad}.  First we will show that
$$\frac{\Lambda}{\Lambda(1-e)\Lambda}\cong\frac{\Tens_{eSe}(eVe)}{\gen{\pi\otimes\pi(R)}}$$
where $\pi:V\onto eVe$ is the obvious map $v\mapsto eve$ of $k$-modules.  This will, in particular, imply that $\Lambda/\Lambda(1-e)\Lambda$ is quadratic.  Our degree $0$ and degree $1$ parts of the isomorphism, $S/S(1-e)S\cong eSe$ and $V/\left(V(1-e)\oplus(1-e)V\right)\cong eVe$, are clear.  We now define a map from $\Lambda$ to the tensor algebra on the right hand side by $v\mapsto eve$ for $v\in V$.  This map is clearly surjective and it is well-defined: if $\sum v_iw_i\in R$ then $\sum ev_ieew_ie\in \pi\otimes\pi(R)$.  As $1=e+(1-e)$ and $v=1v1$ the kernel of this map is the ideal generated by $(1-e)$ and so we have our isomorphism.

Now, as $\Lambda^!=\Tens_S(V^*)/(R^\perp)$, we know that
$$\left(\frac{\Lambda^!}{\Lambda^!(1-e)\Lambda^!}\right)^!\cong \frac{\Tens_{eSe}((eV^*e)^*)}{\gen{(\pi\otimes\pi(R^\perp))^\perp}}$$
and so we need to show that, under our assumptions, the right hand side is isomorphic to $e\Lambda e$.  We have equality in degrees $0$ and $1$ and so, as both algebras are quadratic, it suffices to show that we also have an isomorphism in degree $2$, i.e.,
$$\frac{(eV^*e)^*\otimes_{eSe}(eV^*e)^*}{(\pi\otimes\pi(R^\perp))^\perp}\cong e\left(\frac{V\otimes_SV}{\gen{R}}\right)e\fs$$
We are assuming that $e\Lambda e$ is generated in degree $1$, so the multiplication $(e\Lambda e)_1\otimes_{eSe}(e\Lambda e)_1\to (e\Lambda e)_2$, i.e., the 
composition
$$eVe\otimes_{eSe}eVe\into eV\otimes_S Ve\onto \frac{eV\otimes_S Ve}{\gen{eRe}}\cma$$
is surjective.  Define a map $\psi:(eV^*e)^*\otimes_{eSe}(eV^*e)^*\to (e\Lambda e)_2$ by precomposing with $\psi\otimes\psi$, where $\psi:(e\Lambda e)_2\to eVe$ is induced by the inverse of the evaluation isomorphism $V\arr{\sim}V^{**}$.  Then we calculate the kernel of $\psi$: this is $(eV^*e)^*\otimes_{eSe}(eV^*e)^*\cap(eRe)^**$.  Consider $(\ker \psi)^\perp$: after identifying $V^{**}$ and $V$ we see this is just $\pi\otimes\pi(R^\perp)$, and so $\ker \psi\cong(\pi\otimes\pi(R^\perp))^\perp$ and we have our isomorphism.
\end{pf}
\end{lem}
Note that it is not automatically true that $e\Lambda e$ is generated in degree $1$ or that it is quadratic. For counterexamples consider $\Gamma_3$ and the idempotents $e=e_1$ and $e=e_1+e_2$ respectively, and, using the isomorphism $e\Lambda e\cong\End_\Lambda(\Lambda e)^\op$, apply Lemma \ref{gammaend}.

Now let $$\pi:\Lambda^!\onto \frac{\Lambda^!}{\Lambda^!(1-e)\Lambda^!}$$ be the quotient map and let $M'$ be a $\frac{\Lambda^!}{\Lambda^!(1-e)\Lambda^!}\da\frac{\Lambda^!}{\Lambda^!(1-e)\Lambda^!}$-bimodule.  Using $\pi$ we can inflate $M'$ on both sides to obtain a $\Lambda^!\da\Lambda^!$-bimodule $M$.

\begin{prop}\label{Qinflate}
Suppose $e\Lambda e$ is generated in degree $1$ and quadratic. 
We have two functors
$$Q:\Lambda^!\grmodgr \Lambda^!\to\lin(\Lambda\grprojj \Lambda)$$
 and 
$$Q':\frac{\Lambda^!}{\Lambda^!(1-e)\Lambda^!}\grmodgr\frac{\Lambda^!}{\Lambda^!(1-e)\Lambda^!}\to\lin(e\Lambda e\grprojj e \Lambda e)\fs$$
 Then for $M\in\Lambda^!\grmodgr \Lambda^!$ the inflation of $M'$, 
$$Q(M)=\Lambda e\otimes_{e\Lambda e} Q(M') \otimes_{e\Lambda e} e\Lambda\fs$$

\begin{pf}
Note that $M\cong eMe$ as $\Lambda^!\da\Lambda^!$-bimodules, so in particular, $(M_i)^*\cong e(M_i)^*e$ for each $i\in\Z$.  So in degree $i$
$$Q(M)_i=\Lambda\otimes_S (M_i)^*\otimes_S\Lambda\cong\Lambda\otimes_S e(M_i)^*e\otimes_S\Lambda\cong\Lambda e\otimes_S (M_i)^*\otimes_S e\Lambda\cong\Lambda e\otimes_{eSe} (M_i')^*\otimes_{eSe} e\Lambda$$
and
$$\Lambda e\otimes_{e\Lambda e} Q(M')_i \otimes_{e\Lambda e} e\Lambda=\Lambda e\otimes_{e\Lambda e} e\Lambda e\otimes_{eSe}(M_i')^* \otimes_{eSe}  e\Lambda e \otimes_{e\Lambda e} e\Lambda\cong\Lambda e\otimes_{eSe} (M_i')^*\otimes_{eSe} e\Lambda\cma$$
and it is obvious that the differentials agree.
\end{pf}
\end{prop}

\subsection{Preprojective algebras of type $A$}
Recall that the quadratic dual of $\Gamma_n$ is the preprojective algebra $\Pi_n$ of type $A_n$.  It has an explicit description as follows: it is the quotient of the path algebra of the quiver
$$\xymatrix@=10pt{
  Q^*_n&=&1\ar@/^/[rr]^{x_1} & & 2\ar@/^/[ll]^{y_2}\ar@/^/[rr]^{x_2} & & {}\cdots\ar@/^/[ll]^{y_3}\ar@/^/[rr]^{x_{n-1}} & & n\ar@/^/[ll]^{y_n}
}$$
by the ideal generated by $x_1y_2$, $x_iy_{i+1}-y_ix_{i-1}$ for $2\leq i \leq n-1$, and $y_nx_{n-1}$.  It is graded by path length, and we have identified $x_i$ with $\beta_{i+1}^*$ and $y_j$ with $(-1)^{j-1}\alpha_{j-1}^*$ in the quiver $Q_n$ of $\Gamma_n$.

Consider the algebra surjections
$$\pi^\ell:\Pi_n\onto \frac{\Pi_n}{\Pi_ne_n\Pi_n}\cong \Pi_{n-1}$$
and
$$\pi^r:\Pi_n\onto \frac{\Pi_n}{\Pi_ne_1\Pi_n}\cong \Pi_{n-1}$$
with notation chosen to represent whether the left or right idempotents remain nonzero.
These give us two ways to inflate a $\Pi_{n-1}$-module to a $\Pi_n$-module.

Let $\vec{A}_n$ denote the quiver
$$\vec{A}_n\:\:\:=\:\:\:
1\to 2\to \cdots \to n$$
of Dynkin type $A_n$ with all arrows oriented $i\to i+1$.  Then we have an algebra surjection $\Pi_n\onto k\vec{A}_n$ defined by quotienting out by the ideal generated by all arrows $y_i$.
This gives $k\vec{A}_n$ the structure of a $\Pi_n\da\Pi_n$-bimodule.
\begin{lem}\label{prep-ses}
There is a short exact sequence of graded $\Pi_n\da\Pi_n$-bimodules
$$0\to{}_{\pi^r}(\Pi_{n-1})_{\pi^\ell}\grsh{-1}\into \Pi_n\onto k\vec{A}_n\to0$$
where ${}_{\pi^r}(\Pi_{n-1})_{\pi^\ell}$ denotes the inflation of $\Pi_{n-1}\in\Pi_{n-1}\grmodgr\Pi_{n-1}$ using $\pi^r$ on the left and $\pi^\ell$ on the right.

\begin{pf}
The map $\Pi_n\onto k\vec{A}_n$ is given by the algebra surjection described above and its kernel is the submodule $\gen{y_2,\ldots,y_n}$ of $\Pi_n$, which is generated in degree $1$.  It is easy to check that $e_i\mapsto y_{i+1}$ defines a map ${}_{\pi^r}(\Pi_{n-1})_{\pi^\ell}\grsh{-1}\to\gen{y_2,\ldots,y_n}$ of right modules, and that it respects the left module structure.  Similarly, we have an inverse bimodule map $y_j\mapsto e_{j-1}$, so we are done.
\end{pf}
\end{lem}

The quadratic dual algebras $\Gamma_n$ and $\Pi_n$ (for $n\geq3$) were studied by Brenner, Butler, and King using their theory of almost Koszul duality \cite{bbk}.
We will need the following result:
\begin{prop}[Brenner-Butler-King]\label{preproj-frob}
The preprojective algebra $\Pi_n$ is Frobenius of Gorenstein parameter $n-1$ with Nakayama automorphism $\tau_n^!$.

\begin{pf}
We have an isomorphism of $\Pi_n\da\Pi_n$-bimodules $$\Pi_n\arr{\sim}{}_{\tau_n^!}(\Pi_n)^*$$
by Corollary 4.7 of \cite{bbk}, using the fact that $\tau_n^!$ is its own inverse, and it is easy to calculate the necessary grading shift.
\end{pf}
\end{prop}
Note that the Gorenstein parameter $n-1$ corresponds to the fact that $\Pi_{n}$ is $(n-1,2)$-Koszul \cite[Corollary 4.3]{bbk}.

We will also need a strengthening of Theorem \ref{bbkres}, which is a combination of Theorem 3.15 and Proposition 5.1 of \cite{bbk}.
\begin{thm}[Brenner-Butler-King]\label{bbkres2}
The algebra $\Gamma_n$ is twisted periodic with period $n$, automorphism $\tau_n$, and truncated resolution $Y_n=Q(\Pi_n)$.
\end{thm}

\subsection{Explicit isomorphisms}
Theorem \ref{rzthm} told us how longest elements of $B_{n+1}$ act on the derived categories of the algebras $\Gamma_n$.  The proof of this theorem proceeds by calculating the action of $t_{n+1}$ on the indecomposable projective $\Gamma_n$-modules, showing that two-sided tilting complexes with isomorphic restrictions to one side can only differ by a twist \cite[Proposition 2.3]{rz}, and explicitly determining the outer automorphism groups of the algebras $\Gamma_n$ \cite[Proposition 4.4 and Remark 3]{rz}.

In the rest of this section we present an alternative, explicit proof of Theorem \ref{rzthm}.  Our approach involves working directly with two-sided tilting complexes and uses the theory of almost Koszul duality due to Brenner, Butler, and King \cite{bbk}.  To simplify our notation, from now on, $A$ will denote the algebra $\Gamma_n$.

\begin{lem}\label{lem-gh}
The composition $F_mF_{m-1}\ldots F_2F_1$ is given by tensoring with the complex $$H_m=\cone(G_m\arr{\ev}A)$$ in $\Ch(A\da A)\cma$ where $G_m$ is the complex
$$G_m=\ldots\to0\to
G_{m,m-1}\arr{}\ldots\arr{}G_2\arr{}G_{m,1}\arr{}G_{m,0}\to0\to\ldots$$
 of $A\da A$ bimodules
where $G_{m,i}$ is in degree $i$ and $G_{m,i}=\bigoplus_{j=1}^{m-i}P_{i+j,j}$.  Here,
$P_{i,j}=P_i\otimes_kP_i^\vee$, where $P_i=Ae_i$ is the $i$th indecomposable projective, so the evaluation map is $\ev:\bigoplus_{j=1}^m P_{j,j}\arr{}A$.
The differential is, up to sign, induced by the maps $e_jA\to e_{j+1}A$ and $Ae_i\to Ae_{i-1}$ given by left multiplication by $\beta_{j+1}$ and right multiplication by $\beta_i$, using $P_j^\vee\cong e_jA$.

For example,
$$
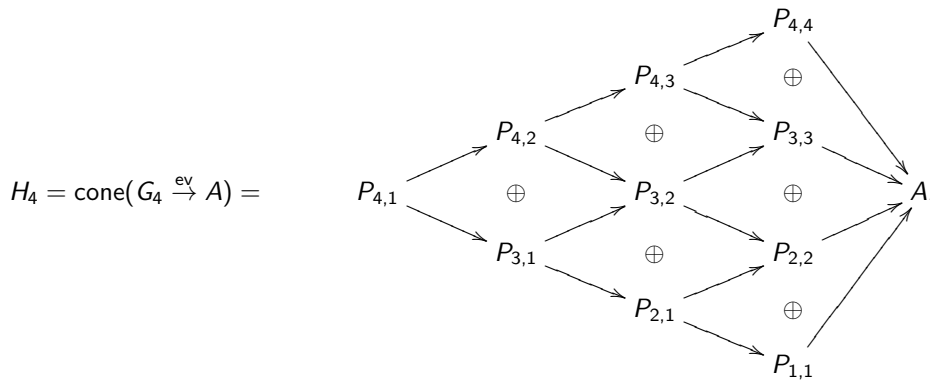
\xymatrix @R=6pt @C=30pt{
  & & & &P_{4,4}\ar[rddd] & \\
  & & &P_{4,3}\ar[ru]\ar[rd] &\oplus & \\
  & &P_{4,2}\ar[ru]\ar[rd] &\oplus &P_{3,3}\ar[rd] & \\
H_4=\cone(G_4\arr{\ev}A)= &P_{4,1}\ar[ru]\ar[rd] &\oplus &P_{3,2}\ar[ru]\ar[rd] &\oplus &A\fs \\
  & &P_{3,1}\ar[ru]\ar[rd] &\oplus &P_{2,2}\ar[ru] & \\
  & & &P_{2,1}\ar[ru]\ar[rd] &\oplus & \\
  & & & &P_{1,1}\ar[ruuu] & 
}$$

\begin{pf}
As $A$ is symmetric, $P_j^\vee\cong P_j^*\cong e_jA^*\cong e_jA$.

We proceed by induction.  Recall that $F_i$ is given by tensoring with the complex $P_{i,i}\arr{\ev}A$ of $A\da A$-bimodules, so the case $m=1$ is clear.

Now suppose that we have shown the statement for $m-1$.  
For $1\leq i,j<m$, $P_{m,m}\otimes_AP_{i,j}$ is isomorphic to $P_m\otimes_k\gen{\beta_m}\otimes_kP_j^\vee\cong P_{m,j}$ if $i=m-1$, and to the zero module otherwise.  So $P_{m,m}\otimes_AG(m-1)$ is the complex
$$P_{m,1}\arr{}P_{m,2}\arr{}\ldots\arr{}P_{m,m-1}\arr{}P_{m,m}$$
with maps as we expect, and $A\otimes_AG(m-1)$ is just $G(m-1)$.  As the differential in $F_i$ is the evaluation map, the isomorphism $P_{m,m}\otimes_AP_{m-1,j}\cong P_m\otimes_k\gen{\beta_m}\otimes_kP_j^\vee$ tells us that these complexes glue together to give the complex $G(m+1)$.
\end{pf}
\end{lem}

Recall that a two-sided tilting complex is a chain complex $X\in\Db(A\da B)$ such that $X\dert_B-:\Db(B)\to\Db(A)$ is an equivalence of triangulated categories \cite{ric-two}.  All the autoequivalences of derived categories we consider in this article are given by tensoring with two-sided tilting complexes, so we can work inside the \emph{derived Picard group}, which we denote $\DPic(A)$.  This is the group of isomorphism classes of two-sided tilting complexes in $\Db(A\da A)$, with group product given by taking the tensor product over $A$ \cite[Definition 3.1]{rz}. 

Our braid group action defines a group morphism
$$\psi_n:B_{n+1}\to\DPic(A)$$
$$\sigma_i\mapsto X_i\fs$$

\emph{Proof of Theorem \ref{rzthm}:}
Recall that $t_{n+1}$ is the positive lift of the longest element in $B_{n+1}$.  We want to determine the action of $t_{n+1}$ on $\Db(A)$, where $A=\Gamma_n$, for all $n\geq 1$.  We check this directly for $n=1$: this is easy.  With a little more work we can also check this directly for $n=2$ and $n=3$.  Then we proceed by induction: assume that $t_n$ acts on $\Db(\Gamma_{n-1})$ as $-_{\tau_{n-1}}[n-1]$ and we will show the corresponding statement for $\Db(A)$.

We want to show that two functors are naturally isomorphic: the shift and twist $-_{\tau_n}[n]$ and the image of $t_{n+1}$ in the group morphism $\varphi_n:B_{n+1}\to \DPic(A)$.  The first is naturally isomorphic to $A_{\tau_n}[n]\otimes_A-$, and the second to some complex of $A\da A$-bimodules obtained by tensoring together complexes $X_i$.  So it is enough to show that these bimodule complexes are isomorphic in $\Db(A\da A)$.

Let $t_n'$ denote the image of $t_n$ in the group monomorphism $B_n\into B_{n+1}$ which sends $\sigma_i\mapsto\sigma_i$.  Let $T_{n+1}$ and $T_n'$ denote the image of $t_{n+1}$ and $t_n'$ in $\psi_n:B_{n+1}\to\DPic(A)$. 
Then we want to show that $T_{n+1}\cong A_{\tau_n}[n]$.

By the inductive description of $t_{n+1}$, we need to show that
$$T_n'X_nX_{n-1}\ldots X_2X_1\cong A_{\tau_n}[n]\fs$$
By our inductive hypothesis we know that $\varphi_{n-1}(t_n)$ acts as a shift and twist on $\Db(\Gamma_{n-1})$, so by the lifting theorem \ref{lifting} $\varphi_{n}(t'_n)$ is given by a periodic twist: we have a distinguished triangle
$$PY_{n-1}P^\vee\arr{}A\to T_n'\rsar{}$$
where $P=P_1\oplus\cdots\oplus P_{n-1}\in A\mMod$, so by Lemma \ref{gammaend} we have $E=\End_A(P)^\op\cong\Gamma_{n-1}$ and $Y_{n-1}$ is the truncated resolution of $\Gamma_{n-1}$ by Theorem \ref{bbkres2}.

By Lemma \ref{lem-gh} we know that $X_nX_{n-1}\ldots X_2X_1\cong H_n$ and we have a distinguished triangle
$$G_n\arr{}A\to H_n\rsar{}$$
so we want to show that
$$T_n'H_n\cong A_{\tau_n}[n]$$
which, as $T_n'$ is a two-sided tilting complex, is equivalent to showing
$$H_n\cong ((T_n')^*)_{\tau_n}[n]\fs$$
We will show that there is a map between these two complexes in the derived category such that the cone is isomorphic to zero.

By Proposition \ref{preproj-frob} we have an isomorphism $\Pi_{n-1}\cong {}_{\tau_{n-1}^!}(\Pi_{n-1})^*\grsh{2-n}$ of $\Pi_{n-1}\da\Pi_{n-1}$-bimodules, so, using that $\tau_{n-1}^!\circ \pi^r=\pi^\ell\circ\tau_n^!$, we have an isomorphism
$${}_{\pi^r}(\Pi_{n-1})_{\pi^\ell}\cong {}_{\tau_n^!}({}_{\pi^l}(\Pi_{n-1}^*)_{\pi^\ell})\grsh{2-n}$$
of $\Pi_n\da\Pi_n$-bimodules.  Applying Proposition \ref{Qprop} \ref{Qtwist} we see that
$$Q({}_{\pi^r}(\Pi_{n-1})_{\pi^\ell})\cong Q({}_{\tau_n^!}({}_{\pi^l}(\Pi_{n-1}^*)_{\pi^\ell})\grsh{2-n})
\cong Q(({}_{\pi^l}(\Pi_{n-1}^*)_{\pi^\ell})\grsh{2-n})_{\tau_n}$$
and Proposition \ref{Qinflate}, which is valid as $n-1>2$, and parts \ref{Qshift} and \ref{Qdual} of Proposition \ref{Qprop} say that this is isomorphic to
$$PQ'(\Pi_{n-1}^*\grsh{2-n})P^\vee_{\tau_n}\cong PQ'(\Pi_{n-1}^*)P^\vee_{\tau_n}[n-2]\grsh{2-n}\cong PY_{n-1}^*P^\vee_{\tau_n}[n-2]\grsh{-n-2}$$
as $\Gamma_{n-1}$ is symmetric of Gorenstein parameter $2$ because it is $(2,n-2)$-Koszul \cite[Proposition 3.11 and Corollary 4.3]{bbk}.

Combining Lemmas \ref{Qexact} and \ref{prep-ses} tells us that we have a short exact sequence
$$0\to Q(k\vec{A}_n)\into Q(\Pi_n)\onto Q({}_{\pi^r}(\Pi_{n-1})_{\pi^\ell}\grsh{-1}) \to0$$
in $\lin(\Gamma_n\grprojj\Gamma_n)$.  We know by Theorem \ref{bbkres2} that $Q(\Pi_n)=Y_n$ and it is easy to see that $Q(k\vec{A}_n)\cong G_n$, so our short exact sequence is
$$0\to G_n\into Y_n\onto  (PY_{n-1}^*P^\vee)_{\tau_n}[n-1]\grsh{-n-1}\to0\fs$$
Forgetting the grading, we get the short exact sequence
$$0\to G_n\into Y_n\onto (PY_{n-1}^*P^\vee)_{\tau_n}[n-1]\to0$$
of chain complexes of $A\da A$-bimodules.

We build a diagram
$$\xymatrix{
 &G_n\ar@{->>}[r]\ar@{^(->}[d] &A\ar@{=}[d]\\
A_{\tau_n}[n-1]\ar@{^(->}[r]\ar@{=}[d] &Y_n\ar@{->>}[r]\ar@{->>}[d] &A\\
A_{\tau_n}[n-1]\ar[r] &(PY_{n-1}^*P^\vee)_{\tau_n}[n-1] &
}$$
which has exact columns: the only non-trivial column comes from the short exact sequence above.  The map in the top row comes from Lemma \ref{lem-gh} and the maps in the middle row come from Theorem \ref{bbkres2}, so it is clear that the top right square commutes.

We now describe the construction of the map $A_{\tau_n}[n-1]\to(PY_{n-1}^*P^\vee)_{\tau_n}[n-1]$.  First, we use the usual periodic twist construction to define a map $PY_{n-1}P^\vee\to A$ from the map $Y_{n-1}\onto E$ given by Theorem \ref{bbkres2}; recall that $E=\End_A(P)^\op\cong\Gamma_{n-1}$.  Then take the dual of this map and pre- and post-compose maps as follows:
$$A\arr{\sim}A^*\to(PY_{n-1}P^\vee)^*\arr{\sim}PY_{n-1}^*P^\vee\fs$$
Finally twist on the right by $\tau_n$ and apply the shift functor $[n-1]$.

Assume for a moment that our diagram commutes.  Then the rows give a short exact sequence
$$0\to H_n\into U\onto ((T_n')^*)_{\tau_n}[n+1]\to 0$$
of chain complexes of $A\da A$-bimodules where $U$ is the complex
$$A_{\tau_n}[n-1]\into Y_n\onto A$$
and so has zero homology.  Therefore we have a distinguished triangle
$$H_n\to 0\to ((T_n')^*)_{\tau_n}[n+1]\rsa$$
and so the map $((T_n')^*)_{\tau_n}[n]\arr{\sim}H_n$ in a rotation of this triangule must be an isomorphism.

Now we show the bottom left square commutes.  Note that $A_{\tau_n}[n-1]\to(PY_{n-1}^*P^\vee)_{\tau_n}[n-1]$ is constructed from the dual of $Y_{n-1}\onto \Gamma_{n-1}$ which is the start of a bimodule resolution and so given by evaluation maps $\bigoplus{P'_{i,i}}\arr{\ev}\Gamma_{n-1}$.  By \cite[Lemma 4.3]{gra1} ($\Delta'\cong\Delta^*$) we can, up to shift and twist, identify $A_{\tau_n}[n-1]\into Y_n$ with the dual of $Y_n\onto A$, which is also given by evaluation maps $\bigoplus{P_{j,j}}\arr{\ev}\Gamma_{n}$.  So as $PP'_{i,i}P^\vee\cong P_{i,i}$ and $Y_n\onto (PY_{n-1}^*P^\vee)_{\tau_n}[n-1]$ is surjective the diagram commutes.
\endpf


\begin{thebibliography}{MMM}

\bibitem[BGS]{bgs} A. Beilinson, V. Ginzburg, and W. Soergel, \emph{Koszul duality patterns in representation theory}, J. Amer. Math. Soc. \textbf{9} (1996), no. 2, 473-527
\bibitem[BoKa]{bk} A. Bondal and M. Kapranov, \emph{Enhanced Triangulated Categories}, Mat. Sb. \textbf{181} (1990), no. 5, 669-683, translation in Math. USSR-Sb. 70 no. 1, 93107
\bibitem[BG]{bg} A. Braverman and D. Gaitsgory, \emph{Poincar\'e-Birkhoff-Witt Theorem for Quadratic Algebras of Koszul Type}, J. Algebra \textbf{181} (1996) 315-328
\bibitem[BuKi]{bking} M. Butler and A. King, \emph{Minimal resolutions of algebras}, J. Algebra \textbf{212} (1999) 323-362
\bibitem[BBK]{bbk} S. Brenner, M. Butler, and A. King, \emph{Periodic algebras which are almost Koszul}, Algebr. Represent. Theory \textbf{5} (2002),  no. 4, 331-367
\bibitem[ERZ]{erz} S. Eilenberg, A. Rosenberg, and D. Zelinsky, \emph{On the dimension of modules and algebras, VIII: dimension of tensor products}, Nagoya Math. J. \textbf{12} (1957), 71-93.
\bibitem[Gra1]{gra1} J. Grant, \emph{Derived autoequivalences from periodic algebras}, Proc. London Math. Soc. (2) \textbf{106} (2013), no. 2, 375-409
\bibitem[Gra2]{gra2} J. Grant, \emph{Derived autoequivalences and braid relations}, Proceedings of the 44th Symposium on Ring Theory and Representation Theory, 50-54, Symp. Ring Theory Represent. Theory Organ. Comm., Okayama, 2012.
\bibitem[Hap]{ha} D. Happel, \emph{On the derived category of a finite-dimensional algebra},
Comment. Math. Helv. 62 (1987), no. 3, 339-389
\bibitem[HK]{hk} R. S. Huerfano and M. Khovanov, \emph{A category for the adjoint representation}, J. Algebra \textbf{246} (2001), no. 2, 514-542
\bibitem[Ive]{ive} B. Iversen \emph{Octahedra and braids}, Bulletin de la Soci\'et\'e Math\'ematique de France \textbf{114} (1986) 197-213
\bibitem[KN]{kn} B. Keller and A. Neeman, \emph{The connections between May's axioms for a triangulated tensor product and Happel's description of the derived category of the quiver $D_4$}, Doc. Math. J. DMV \textbf{7} (2002), 535-560
\bibitem[Mat]{matsumoto} H. Matsumoto, \emph{G\'en\'erateurs et relations des groupes de Weyl g\'en\'eralis\'es}, C. R. Acad. Sci. Paris \textbf{258} (1964), 3419-3422
\bibitem[May]{m} J. P. May, \emph{The additivity of traces in triangulated categories}, Adv. Math. \textbf{163} (2001), no. 1, 34-73
\bibitem[MOS]{mos} V. Mazorchuk, S Ovsienko, and C. Stroppel, \emph{Quadratic duals, Koszul dual functors, and applications}, Trans. Amer. Math. Soc. \textbf{361} (2009), no. 3, 1129-1172. 
\bibitem[Pri]{pri} S. Priddy, \emph{Koszul resolutions}, Trans. Amer. Math. Soc. \textbf{152} (1970), 39-60
\bibitem[Ric]{ric-two} J. Rickard, \emph{Derived equivalences as derived functors}, J. London Math. Soc. (2) \textbf{43} (1991) 37-48
\bibitem[RZ]{rz} R. Rouquier and A. Zimmermann, \emph{Picard groups for derived module categories}, Proc. London Math. Soc. (3) \textbf{87} (2003), no. 1, 197-225
\bibitem[Sei]{sei} P. Seidel, \emph{Lagrangian homology spheres in $(A_m)$ Milnor fibres}, 	arXiv:1202.1955v2 [math.SG], preprint 2012
\bibitem[ST]{st} P. Seidel and R. Thomas, \emph{Braid group actions on derived categories of coherent sheaves}, Duke Math. J. \textbf{108} (2001), no. 1, 37-108
\end{thebibliography}
\end{document}